\documentclass[12pt,letterpaper]{amsart}

\allowdisplaybreaks
\usepackage{mathpazo}
\usepackage{amsmath,amsfonts,amssymb}
\usepackage{amsrefs}
\usepackage{amsthm}
\usepackage{latexsym,amsmath,amssymb,amsfonts}
\usepackage{rotating}
\usepackage{mathrsfs}
\usepackage{xypic} \xyoption{all}
\usepackage{amscd}
\usepackage{hyperref} 
\usepackage{euscript}
\usepackage{hhline}
\usepackage{graphicx,epstopdf}
\usepackage{epsfig}
\usepackage{xcolor}
\usepackage{textcomp}
\usepackage[all,color]{xy}
\usepackage{spectralsequences}
\usepackage{sseq}
\usepackage{tikz}

\newlength{\fighskip} \fighskip=2pt
\newlength{\figvskip} \figvskip=3pt

\usepackage{hyperref}

\advance\textheight by \topskip
\numberwithin{equation}{section}

\newcommand{\C}{\mathbb{C}}

\newcommand{\W}{\mathcal{W}}

\DeclareMathOperator{\End}{End}

\DeclareMathOperator{\Sym}{Sym}
\DeclareMathOperator{\Hom}{Hom}

\newcommand{\A}{\mathcal A}

\renewcommand{\W}{\mathcal W}

\theoremstyle{plain}
\newtheorem{thm}{Theorem}[section]
\newtheorem{thm-defn}{Theorem/Definition}[section]
\newtheorem{lem}[thm]{Lemma}
\newtheorem{lem-defn}[thm]{Lemma/Definition}
\newtheorem{prop}[thm]{Proposition}
\newtheorem{cor}[thm]{Corollary}

\theoremstyle{definition}
\newtheorem{defn}[thm]{Definition}

\newtheorem{eg}[thm]{Example}

\theoremstyle{remark}
\newtheorem{rmk}[thm]{Remark}

\allowdisplaybreaks[4]  

\begin{document}
\title[Quantization of K\"ahler manifolds via differential operators]{Quantization of K\"ahler manifolds via differential operators}

\author[Chan]{Kwokwai Chan}
\address{Department of Mathematics\\ The Chinese University of Hong Kong\\ Shatin \\ Hong Kong}
\email{kwchan@math.cuhk.edu.hk}

\author[Leung]{Naichung Conan Leung}
\address{The Institute of Mathematical Sciences and Department of Mathematics\\ The Chinese University of Hong Kong\\ Shatin \\ Hong Kong}
\email{leung@math.cuhk.edu.hk}

\author[Li]{Qin Li}
\address{Shenzhen Institute for Quantum Science and Engineering, Southern University of Science and Technology, Shenzhen, China}
\email{liqin@sustech.edu.cn}

\author[Yau]{Yutung Yau}
\address{Kavli Institute for the Physics and Mathematics of the Universe (WPI), The University of Tokyo Institutes for Advanced Study, The University of Tokyo, Kashiwa, Chiba 277-8583, Japan}
\email{yu-tung.yau@ipmu.jp}


\begin{abstract}
In this paper, we study the quantization of classical observables (i.e., functions) on a K\"ahler manifold $X$ as differential operators acting on holomorphic sections of tensor powers $L^{\otimes k}$ of the pre-quantum line bundle $L$. We prove two global results as follows.
\begin{itemize}
    \item For a general smooth function $f \in C^\infty(X)$, we construct higher order generalizations of Kostant--Souriau's pre-quantum differential operators using our Fedosov-type constructions of Bargmann--Fock sheaves in \cite{ChaLeuLi2022b, ChaLeuLi2021, ChaLeuLi2023}. We prove that these differential operators are asymptotic to the Berezin--Toeplitz operators $T_{f,k}$ acting on the Hilbert space $\mathcal{H}_k = H^0(X, L^{\otimes k})$ as $k \to \infty$.
    \item If a smooth function $f \in C^\infty(X)$ is furthermore the symbol of a level $k$ quantizable function (as defined in \cite{ChaLeuLi2023}), then we prove that the associated Berezin--Toeplitz operators $T_{f,k}$ are all holomorphic differential operators. Conversely, Berezin--Toeplitz operators that are holomorphic differential operators all arise in this way. This gives a complete characterization of when Berezin--Toeplitz operators are holomorphic differential operators.
\end{itemize}
To prove these results, we establish new orthogonality relations which generalize the classical Tuynman's Lemma, 
and employ various differential-geometric and analytic technqiues such as H\"ormander's estimates.


\end{abstract}

\maketitle


\section{Introduction}

\subsection*{Background}

The quantization of a classical mechanical system turns the phase space $X$ (a symplectic manifold) into a Hilbert space $\mathcal{H}$ and classical observables (i.e., functions on $X$) into operators on $\mathcal{H}$. 

Consider the example of a flat space $X = \mathbb{R}^{2n}$. If we regard $X = T^*\mathbb{R}^n$, where $(x_1,\dots,x_n)$ are the base coordinates and $(p_1,\dots,p_n)$ are the fiber coordinates, the {\em canonical quantization} assigns $x_j \mapsto \hat{x}_j = m_{x_j}$ (i.e., multiplication by $x_j$) and $p_j \mapsto \hat{p}_j = -i\hbar\frac{\partial }{\partial x_{j}}$. These operators on the Hilbert space $\mathcal{H}_{\mathbb{R}} = L^2(\mathbb{R}^{n})$ (or the Schwartz space) satisfy the uncertainty principle $[\hat{x}_{j},\hat{p}_{k}] \sim i\hbar\delta _{jk}$ and define a star product $\star_M$ on $\mathcal{H}_\mathbb{R}$ called the {\em Moyal product} via composition of differential operators.

If instead we regard $X = \mathbb{C}^n$ with holomorphic coordinates $(z_1,\dots,z_n)$, the Hilbert space is given by the {\em Bargmann--Fock space} $\mathcal{H}_\mathbb{C}$ of holomorphic functions $s(z)$ that are $L^{2}$ with respect to $\left\Vert s\right\Vert^{2}=\int_{X}\left\vert s(z) \right\vert ^{2}e^{-\left\vert z\right\vert ^{2}/2}$. We then quantize by sending $z_j \mapsto \hat{z}_j = z_j \cdot$ and $\bar{z}_j \mapsto \hat{\bar{z}}_j = \frac{\hbar}{\sqrt{-1}} \frac{\partial}{\partial z_j}$. This defines another star product $\star_W$, via composition of holomorphic differential operators, called the {\em Wick product}.

In the flat case, it is clear that the space of functions on $T^*\mathbb{R}^n$ which are polynomials along the fibers, equipped with a star product $\star$, genuinely acts on the Hilbert space $\mathcal{H}$ in geometric quantization as differential operators. Similarly, functions on $\mathbb{C}^n$ which are polynomials in the anti-holomorphic variables $\bar{z}_1,\cdots, \bar{z}_n$ act on $\mathcal{H}_\mathbb{C}$ as holomorphic differential operators. 
It is natural to ask whether such a quantization via differential operators can be generalized to curved manifolds. In this paper, we will focus on the K\"ahler setting, as in the previous series of works \cite{ChaLeuLi2022b, ChaLeuLi2021, ChaLeuLi2023}.

For a general compact symplectic manifold $(X, \omega)$, the Hilbert space of the quantum mechanical system is given by Kostant--Souriau's geometric quantization \cite{Kostant70, Souriau70}. One needs to choose a {\em pre-quantum line bundle}, namely, a complex line bundle $L$ equipped with a connection $\nabla$ such that $\omega = \sqrt{-1}\cdot F_\nabla$, where $F_\nabla = \nabla^2$ is the curvature. 
One also needs to choose a polarization. For a compact K\"ahler manifold $(X,\omega,J)$, a natural choice is the complex polarization induced by the complex structure $J$. Then, for any positive integer $k$, the {\em level $k$ Hilbert space} is defined as the space 
$$\mathcal{H}_k := H^0(X, L^{\otimes k})$$
of holomorphic sections of the $k$-th tensor power $L^{\otimes k}$ of the pre-quantum line bundle.

To quantize a function $f \in C^\infty(X)$, there are various approaches including the Fedosov quantization \cite{Fed1994} and the Berezin--Toeplitz quantization \cite{Berezin75}. In this paper, we will combine these approaches together to show how smooth functions on a general compact K\"ahler manifold $(X,\omega,J)$ can be quantized as differential operators. 
Ideally, we require that a quantum operator should act on the Hilbert spaces $\mathcal{H}_k$ obtained in geometric quantization. 
For a general smooth function $f\in C^\infty(X)$, although the associated operator only maps $s \in \mathcal{H}_k = H^0(X, L^{\otimes k})$ to a smooth section of $L^{\otimes k}$ in general, we will prove that the outputs are actually {\em asymptotically close} to the subspace of holomorphic sections. 
When $f \in C^\infty(X)$ is in addition (the symbol of) a so-called {\em quantizable function}, we will prove that the associated operator is exactly {\em equal to} a holomorphic differential operator; indeed, the resulting operator acting on $\mathcal{H}_k$ can be identified with the corresponding Berezin--Toeplitz operator $T_{f,k}$. Below is a more detailed description of our main results.




\subsection*{Quantization in terms of higher Kostant--Souriau operators}

Applying the Fedosov-type constructions in \cite{ChaLeuLi2022b, ChaLeuLi2021, ChaLeuLi2023}, we will define differential operators (see Definition \ref{definition: higher-Kostant-Souriau-operator})
$$P_{f, k ,m}: H^0(X, L^{\otimes k})\rightarrow \Gamma(X, L^{\otimes k})\quad \text{for }m\geq 0, k\geq 1,$$
associated to a general smooth function $f \in C^\infty(X)$; here, $\Gamma(X, L^{\otimes k})$ denotes the space of smooth sections of $L^{\otimes k}$.
We call these {\em higher Kostant--Souriau operators} as they are higher order analogues of the classical Kostant--Souriau pre-quantum differential operators.

To construct the operators $P_{f,k,m}$'s, recall that in \cite{ChaLeuLi2022b} we applied Fedosov's geometric construction \cite{Fed1994} to obtain a class of {\em Fedosov flat connections} on the (formal) Weyl bundle $\mathcal{W}_{X,\C}[[\hbar]]$ of a K\"ahler manifold $X$, which we denote by $D_{BT}$ (here the subscript ``BT'' stands for ``Berezin--Toeplitz'').\footnote{In fact, we can obtain all Wick-type deformation quantization using our construction \cite{ChaLeuLi2022b}.}
These Fedosov connections have some distinctive features, e.g., they are compatible with the complex polarization on $X$ and are quantizations of the natural $L_\infty$ structure on a K\"ahler manifold \cite{Kap1999}. Furthermore, as in Fedosov \cite{Fed1994}, the space of $D_{BT}$-flat sections in $\mathcal{W}_{X,\C}[[\hbar]]$ is naturally isomorphic to $C^\infty(X)[[\hbar]]$. In particular, for every smooth function $f\in C^\infty(X)$, we have a corresponding flat section $O_f$ of $\mathcal{W}_{X,\C}[[\hbar]]$.

In \cite{ChaLeuLi2021}, we extended the above to construct Fedosov connections $D_{BF,k}$ (here the subscript ``BF'' stands for ``Bargmann--Fock'') on the {\em Bargmann--Fock sheaf} $\mathcal{F}_k = \mathcal{W}_X\otimes_{\mathcal{O}_X}L^{\otimes k}$. 
We also proved that the space of $D_{BF, k}$-flat sections is canonically isomorphic to the Hilbert space $\mathcal{H}_k = H^0(X, L^{\otimes k})$. In other words, for every holomorphic section $s\in H^0(X, L^{\otimes k})$, there exists a (non-formal) section $O_s$ of $\mathcal{W}_X\otimes_{\mathcal{O}_X}L^{\otimes k}$ flat under $D_{BF,k}$. This produces a module sheaf over the sheaf of deformation quantization algebras.

To have any sensible action of the Weyl bundle $\mathcal{W}_{X,\C}[[\hbar]]$ on the Bargmann--Fock sheaf $\mathcal{F}_k$, we must take the non-formal evaluation $\hbar = \tfrac{\sqrt{-1}}{k}$. In order to do so, we can only allow polynomial dependence on $\hbar$ to avoid divergence issues (as otherwise we may have possibly divergent power series in $\tfrac{1}{k}$).
Actually, if we look more carefully on the flat section $O_s$ of $\mathcal{F}_k$ and the Bargmann--Fock action, we also only allow polynomial dependence on the fiber coordinates $\bar{y}^j$'s of $\mathcal{W}_{X,\mathbb{C}}$ since each such variable acts as $\frac{\sqrt{-1}}{k}\cdot\omega^{\bar{j}i}\cdot\partial_{z^i}$ and thus a power series in $\bar{y}^j$'s may also lead to possibly divergent power series in $\tfrac{1}{k}$. 

This motivates us to define a weight on $\mathcal{W}_{X,\C}$, as the sum of the degrees of $\hbar$ and $\bar{y}^j$'s in a monomial. We call this the {\em polarized weight} as it only counts the degrees of the variables $\bar{y}^j$'s but not that of the variables $y^i$'s. 
In particular, for each flat section $O_f$ of $\mathcal{W}_{X,\C}[[\hbar]]$, we have a weight decomposition
\begin{equation*}\label{equation: formal-sum-Of}
    O_f=\sum_{m \geq 0}(O_f)_m.
\end{equation*}
For each $m\geq 0$, we obtain a well-defined {\em level $k$ Bargmann--Fock action} (here we take the evaluation $\hbar=\tfrac{\sqrt{-1}}{k}$ in $(O_f)_m$)
$$
(O_f)_m\circledast_k O_s. 
$$

\begin{defn}[see Definitions \ref{definition: differential-operators-finite-weight-section-Weyl} and \ref{definition: higher-Kostant-Souriau-operator} in Section \ref{subsection: Kostant-Souriau-operators}]\label{definition: differential-operators-intro}
	Let $f\in C^\infty(X)$ be a smooth function. For any weight $m \geq 0$ and any level $k > 0$, the differential operator
	$$
		P_{f,k,m}: H^0(X, L^{\otimes k})\rightarrow \Gamma(X, L^{\otimes k}),
	$$
	which we call a {\em higher Kostant--Souriau operator}, is defined as the composition of the symbol map with the Bargmann--Fock action
	$$
	P_{f,k,m}(s) := \sigma_{L^{\otimes k}}\left((O_f)_m\circledast_k O_s\right).
	$$
\end{defn}

The operator $P_{f,k,m=0}$ is nothing but $m_f$, the multiplication by $f$. For $m=1$, a simple computation shows that for $s \in H^0(X, L^{\otimes k})$,
$$
P_{f,k,m=1}(s) = \frac{1}{k}\cdot\left(\sqrt{-1} \nabla_{X_f}+\Delta f\right) (s).
$$
Recall that the classical {\em Kostant--Souriau pre-quantum operator} is defined as
\begin{equation}\label{equation: pre-quantum-operator}
	H_{f,k} := \frac{\sqrt{-1}}{k}\cdot\nabla_{X_f}+m_f, 
\end{equation}
where $X_f$ denotes the Hamiltonian vector field associated to $f$ and $\nabla_{X_f}$ is the covariant derivative along $X_f$.
So the operator $P_{f,k, \leq1} := P_{f,k,m=0} + P_{f,k,m=1}$ only differs from the classical Kostant--Souriau pre-quantum operator $H_{f,k}$ by a correction term $\tfrac{1}{k} \Delta f$;
see Example \ref{example: weight-1-differential-operator}.  This explains the name {\em higher Kostant-Souriau operators} for general $P_{f,k,\leq l}$'s. 

For $m \geq 2$, the operators $P_{f,k,m}$ involve higher covariant derivatives for tensor fields (Definition \ref{definition: higher-covariant-derivative}), but they are quite computable. Since the Fedosov construction is an iterative process, it is possible to write down an algorithm (using computer codes) to compute the operators $P_{f,k,m}$. See Example \ref{example: weight-2-differential-operator} for a computation of $P_{f,k,m=2}$ by hand.



The higher Kostant--Souriau operators $P_{f,k,m}$ satisfy two fundamental properties. The first is a set of {\em orthogonality relations}.
\begin{thm}[=Theorem \ref{theorem: orthogonality-relations}]\label{theorem: generalized-Tuyman-lemma}
	For every $m\geq 1$ and any holomorphic section $s\in \mathcal{H}_k = H^0(X, L^{\otimes k})$, the (smooth) section $P_{f,k,m}(s)$ lives in the orthogonal complement $\mathcal{H}_k^\perp\subset\Gamma(X, L^{\otimes k})$, where $\perp$ is with respect to the hermitian metric $h$ on $L^{\otimes k}$.
\end{thm}

Recall that the {\em Berezin--Toeplitz operator} \cite{Berezin75} associated to a smooth function $f \in C^\infty(X)$ is defined as 
$$T_{f,k} := \Pi_k \circ m_{f},$$
where $m_f$ is multiplication by $f$ and $\Pi_k :\Gamma(X, L^{\otimes k})
\rightarrow H^{0}(X, L^{\otimes k})$ is the orthogonal projection.
In \cite{Tuynman87}, Tuynman showed that the Kostant--Souriau pre-quantum operator $H_{f,k}$ is related to the Berezin--Toeplitz operator $T_{f,k}$ by
$\Pi_k \circ H_{f,k} = T_{f - \frac{1}{k} \Delta f}$, 
where $\Delta f$ is the Laplacian of the function $f$ with respect to the K\"ahler metric defined by $\omega$. This is called {\em Tuynman's Lemma}. An equivalent way to state this result is as follows.
\begin{lem}[Tuynman's Lemma \cite{Tuynman87}]\label{lemma: Tuynman}
	For every smooth function $f$ and any two holomorphic sections $s_1,s_2\in H^0(X, L^{\otimes k})$, we have the following orthogonality relation
	$$
	\left\langle \frac{\sqrt{-1}}{k}\nabla_{X_f}s_1, s_2 \right\rangle = \left\langle -\frac{1}{k}\Delta f\cdot s_1, s_2 \right\rangle.
	$$
\end{lem}
\noindent Theorem \ref{theorem: generalized-Tuyman-lemma} is therefore a vast generalization of Tuynman's Lemma \ref{lemma: Tuynman}.

The second key property of the higher Kostant--Souriau operators is the estimates given in the following proposition (see Propositions \ref{proposition: norm-estimated-differential-operators-f} and \ref{proposition: estimate-partial-bar-D-f-action-s}).
\begin{prop}
	We have the following estimates for the higher Kostant--Souriau operators: there exist constants $\tilde{C}_{f,m}$ and $C_{f,m}$, depending only on the function $f$ and the weight $m$, such that 
	\begin{equation}\label{equation: norm-estimate-differential-operators}
		\lVert P_{f,k,m} \rVert_{op} \leq C_{f,m}\cdot k^{-m/2}. 
	\end{equation}	
and 
\begin{equation}\label{equation: estimate-partial-bar-D-f-action-s}
	\lVert \bar{\partial}(P_{f,k,\leq m}s) \rVert \leq \tilde{C}_{f,m}\cdot k^{-\frac{m}{2}}\cdot \lVert s \rVert.
\end{equation}
Here, $\lVert\cdot\rVert$ denotes the $L^2$-norms on $\Gamma(X, L^{\otimes k})$ and $\lVert \cdot \rVert_{op}$ denotes the operator norm.
\end{prop}
Note that the input of the operator $P_{f,k,m}$ is a holomorphic section, but its output is in general only a {\em smooth} section of $L^{\otimes k}$. We need equation \eqref{equation: estimate-partial-bar-D-f-action-s} in Proposition \ref{proposition: estimate-partial-bar-D-f-action-s} to justify the claim that $P_{f,k,m}(s)$ is ``close to being holomorphic''. 
The proof relies essentially on the Fedosov equation for flat sections in the Weyl and Bargmann--Fock bundles, together with some subtle norm estimates of the higher covariant derivatives (see Theorem \ref{theorem: estimate-higher-covariant-derivative}).

Combining these two properties together, we have a more intuitive understanding of the higher Kostant--Souriau differential operators $P_{f,k,m}$'s. 
First of all,
equation \eqref{equation: norm-estimate-differential-operators} implies that the sequence of operators $\{P_{f,k,m}\}_m$ have smaller and smaller norms as $m$ grows.
The first term 
$$P_{f,k,0}(s)=f\cdot s$$
is simply multiplication by $f$. 
For $m\geq 1$, the orthogonality relations in Theorem \ref{theorem: generalized-Tuyman-lemma} say that 
$$P_{f,k,m}\in\mathcal{H}_k^\perp.$$
As a result, the infinite formal sum
\begin{equation} \label{equation: formal-sum-differential-operators}
	 P_{f,k,0}+P_{f,k,1}+P_{f,k,2}+\cdots. 
\end{equation}
can be viewed as giving {\em quantum corrections} to $f\cdot s$ by terms that are orthogonal to $\mathcal{H}_k$, so that its image is {\em asymptotically close} to the Hilbert space $\mathcal{H}_k$ as $k \to \infty$. We conclude that \eqref{equation: formal-sum-differential-operators} should precisely give the {\em orthogonal projection} of $f\cdot s$ in $\mathcal{H}_k$, namely, it is exactly the action of the Berezin--Toeplitz operator $T_{f,k}$ on $s \in H^0(X, L^{\otimes k})$.


This heuristic reasoning is justified by the following comparison between the higher Kostant--Souriau operators and the Berezin--Toeplitz operators, whose detailed proof will be a combination of our Fedosov-type constructions in \cite{ChaLeuLi2022b, ChaLeuLi2021, ChaLeuLi2023} and the H\"ormander's estimates in complex analysis.
\begin{thm}[=Theorem \ref{theorem: main}]\label{theorem: main-result-estimate}
	For any smooth function $f\in C^\infty(X)$ and $m\geq0$, there exists a constant $C_{f,m}$, depending only on $f$ and $m$, such that
	$$
	\left\Vert T_{f,k}-\sum_{l=0}^m P_{f,k,l}\right\Vert_{op} \leq C_{f,m}\cdot\frac{1}{k^{\frac{m+1}{2}}}. 
	$$
	Here $\lVert\cdot\rVert_{op}$ denotes the operator norm.  In other words, the Berezin--Toeplitz operator $T_{f,k}$ is asymptotic to an {\em infinite} formal sum of differential operators, which can be written as 
	$$
	T_{f,k}\sim P_{f,k,0}+P_{f,k,1}+P_{f,k,2}+\cdots. 
	$$
\end{thm}
Guillemin \cite{Guillemin95}, Bordemann--Meinrenken--Schlichenmaier \cite{BorMeiSch1994} and Schlichenmaier \cite{Schlichenmaier00} famously proved that there exist bi-differential operators $C_i(-,-)$ such that the following estimates hold:
\begin{equation}\label{equation: BT-star-product-asymptotic}
	\left\Vert T_{f,k}\circ T_{g,k} - \sum_{i=0}^{N-1}\left(\frac{1}{k}\right)^iT_{C_i(f,g),k}\right\Vert_{op} \leq K_N(f,g)\left(\frac{1}{k}\right)^N;
\end{equation}
here $K_N(f,g)$ are constants that are independent of $k$ and $\Vert \cdot \Vert_{op}$ denotes the operator norm.  This defines the {\em Berezin--Toeplitz star product} $\star_{BT}$ and the {\em Berezin--Toeplitz deformation quantization algebra} $(C^\infty(X)[[\hbar]],\star_{BT})$. 
The Berezin--Toeplitz quantization is one of the most studied quantization schemes and it has numerous applications \cite{Borthwick-Uribe96, Borthwick-Paul-Uribe98, Borthwick-Uribe03, Zelditch98, Bordemann-Waldmann97, Karabegov-Schlichenmaier01, Charles03, ChaPol2017, Ma-Marinescu12, Andersen12, Ioos-Kazhdan-Polterovich21}.

Theorem \ref{theorem: main-result-estimate} can be regarded as a ``module version'' of the estimate \eqref{equation: BT-star-product-asymptotic} in the following sense: equation \eqref{equation: BT-star-product-asymptotic} says that Berezin-Toeplitz operators are not closed under composition in the strict sense, but only asymptotically as $k\rightarrow\infty$. Similarly, outputs of the higher Kostant--Souriau operators are only smooth sections that are holomorphic asymptotically as $k\rightarrow\infty$; see Remark \ref{remark: module-version-of-BMS}.

We may further illustrate the above results by a picture (see Figure \ref{figure: asymptotics}) which describes the relations between the operators. We have the following infinite formal series of differential operators:
$$
P_{f,k,0}+P_{f,k,1}+\cdots,
$$
each term of which has its norm $\lVert P_{f,k,m} \rVert_{op} =  O(k^{-\frac{m}{2}})$. By the orthogonality relations in Theorem \ref{theorem: orthogonality-relations}, we have $P_{f,k,m}(s)\in\mathcal{H}_k^\perp$ for each $m\geq 1$ and any $s\in\mathcal{H}_k$. Modulo convergence issues, terms in the above infinite series are {\em corrections} of $P_{f,k,0} = m_f$ (multiplication by $f$) given by a sequence of operators $\mathcal{H}_k \to \mathcal{H}_k^\perp$ with decreasing operator norms (in the asymptotic sense).

\begin{center}
\begin{figure}
	\begin{tikzpicture}
		\draw[blue] (-5, 0) -- (2, 0);
		\node at (2.5, 0) [right] {\textcolor{blue}{$H^0(X, L^{\otimes k})$}};
		
		\node at (-4, 0) {$\bullet$};
		\node at (-4, -.1) [below] {$s$};

		\node at (0, 5) {$\bullet$};
		\node at (-.5, 5) [left] {$\textcolor{gray}{f \cdot s =} P_{f, k, 0} (s)$};
		\node at (2.5, 5) [right] {$\Gamma(X, L^{\otimes k})$};

		\draw[gray, dashed] (0, 5) -- (0, 0);
		\draw[gray] (-.2, 0) -- (-.2, .2) -- (0, .2);
		\node at (0, 0) {$\bullet$};
		\node at (0, -.1) [below] {$T_{f, k} (s)$};

		\draw[orange] (.2, 4.9) -- (.3, 4.9) -- (.3, 2.6) -- (.2, 2.6);
		\node at (.3, 3.75) [right] {\textcolor{orange}{$P_{f, k, 1}(s)$}};
		\node at (0, 2.5) {$\bullet$};
		\node at (-.5, 2.5) [left] {$\textcolor{gray}{H_{f, k}(s) + \tfrac{1}{k} (\Delta f) \cdot s =} P_{f, k, \leq 1} (s)$};

		\draw[orange] (.2, 2.4) -- (.3, 2.4) -- (.3, 1.4) -- (.2, 1.4);
		\node at (.3, 1.9) [right] {\textcolor{orange}{$P_{f, k, 2}(s)$}};
		\node at (0, 1.5) {$\bullet$};
		\node at (-.5, 1.5) [left] {$P_{f, k, \leq 2} (s)$};

		\node[orange] at (.3, 1.1)  {$\vdots$};

	\end{tikzpicture}
    \caption{Asymptotic relations between operators}
    \label{figure: asymptotics}
\end{figure}
\end{center}

Theorem \ref{theorem: main-result-estimate} is thus saying that, for a general smooth function $f\in C^\infty(X)$, the Berezin--Toeplitz operator $T_{f,k} = \Pi_k\circ m_f$ is asymptotic to the above series of differential operators as $k\rightarrow\infty$. 
Since differential operators are local, an immediate consequence of this theorem is the {\em asymptotic locality} of Berezin--Toeplitz operators. This also implies that $P_{f,k,0}+P_{f,k,1}+\cdots$ asymptotically maps $\mathcal{H}_k$ to itself. Thus the higher Kostant--Souriau operators deserve to be called {\em quantum operators}. 

\subsection*{Quantizable functions}

For a general smooth function $f\in C^\infty(X)$, the classical Kostant--Souriau pre-quantum operator $H_{f,k}$ preserves the Hilbert spaces $\mathcal{H}_k$'s if and only if the Hamiltonian vector field $X_f$ preserves the polarization. Functions satisfying this condition are classically called {\em quantizable functions} or {\em polarization-preserving functions} in the Kostant--Souriau's geometric quantization picture.

This condition is very restrictive (see Remark \ref{remark: quantizable functions}), and quantizable functions obtained this way do not even form an algebra because $H_{f,k}$'s are all just differential operators of order one.
In \cite{ChaLeuLi2023}, we introduced a new notion of {\em quantizable functions}, generalizing the classical one in Kostant--Souriau's geometric quantization. According to the definition (see Definition \ref{definition: quantizable functions}), a {\em quantizable function of level $k$} is a Fedosov flat section $\alpha\in\Gamma(X, \W_X \otimes \Sym \overline{T^*X})$; here $\W_X = \widehat{\Sym} T^*X$ denotes the (holomorphic) Weyl bundle on $X$. 

This new notion of quantizable functions defined in \cite{ChaLeuLi2023} is a vast generalization of the classical one in geometric quantization. For example, the differential operator associated to a quantizable function (in the new sense) can be of arbitrarily high orders. Furthermore, our quantizable functions do form an algebra. 

The second main result of this paper says that the Berezin--Toeplitz operator $T_{f,k}$ acts on the Hilbert space $\mathcal{H}_k$ as a {\em genuine differential operator} precisely when $f$ is the symbol of a level $k$ quantizable function.
\begin{thm}[=Theorem \ref{theorem: Toeplitz-holomorphic-differential-operator}]\label{theorem: quantizable-functions}
	Suppose a smooth function $f \in C^\infty(X)$ is the symbol of a quantizable function of level $k$, i.e., $f = \sigma(\alpha)$ for some section $\alpha\in\Gamma(X, \W_X \otimes \Sym \overline{T^*X})$ flat under the non-formal Fedosov connection $D_{BT,k}$. Then the associated Berezin--Toeplitz operator $T_{f,k}$ is equal to the holomorphic differential operator $P_{\alpha,k}$ on $\mathcal{H}_k$ associated to $\alpha$ (Definition \ref{definition: differential-operators-finite-weight-section-Weyl}).
	
	Conversely, if a Berezin--Toeplitz operator is a holomorphic differential operator, then it can be written as $T_{f,k}$ for some $f$ that is the symbol of a level $k$ quantizable function. 
\end{thm}


The holomorphic differential operator $P_{\alpha,k}$ associated to quantizable function $\alpha$ can locally be described as follows (see Definition \ref{definition: differential-operators-finite-weight-section-Weyl} and Section \ref{subsection: Kostant-Souriau-operators} for details).
When the smooth function $f \in C^\infty(X)$ is of the form $\sigma(\alpha)$ for some level $k$ quantizable function $\alpha$, $f$ is {\em locally} the evaluation of a {\em formal} quantizable function (Definition \ref{definition: formal-quantizable-function}) $\tilde{f} \in C^\infty({X})[[\hbar]]$ at $\hbar = \tfrac{\sqrt{-1}}{k}$ (so that $\alpha$ is locally the same evaluation of $O_{\tilde{f}}$). Then {\em locally} $P_{\alpha,k}$ can be written as a {\em finite} sum of differential operators
\begin{equation}\label{equation: finite-sum-differential-operators}
	P_{\tilde{f},k,0}+P_{\tilde{f},k,1}+\cdots+P_{\tilde{f},k,N},
\end{equation}
where $N$ is the weight of $\tilde{f}$ as a formal quantizable function. This is in sharp contrast with Theorem \ref{theorem: main-result-estimate} where we need an infinite formal sum.

Modulo technical details, the proof of Theorem \ref{theorem: quantizable-functions} goes roughly as follows:
On the one hand, for a quantizable function $f$ and $s\in\mathcal{H}_k$, the finite sum $f\cdot s+P_{f,k,1}(s)+P_{f,k,2}(s)+\cdots$ remains a holomorphic section of $L^{\otimes k}$. This follows from the compatibility between Fedosov connections on the Weyl bundle and Bargmann--Fock sheaf \cite{ChaLeuLi2021}. 
On the other hand, since $P_{f,k,1}(s)+P_{f,k,2}(s)+\cdots\in\mathcal{H}_k^{\perp}$, it follows that $f\cdot s+P_{f,k,1}(s)+P_{f,k,2}(s)+\cdots$ is exactly equal to the orthogonal projection of $f\cdot s$ to $\mathcal{H}_k$, i.e., $T_{f,k}(s)$. 

\begin{rmk}
	Theorems \ref{theorem: main-result-estimate} and \ref{theorem: quantizable-functions} are closely related to very recent work of Andersen \cite{Andersen24}, where he announced that analytic functions on a compact K\"ahler manifold can be quantized to operators on the Hilbert spaces $\mathcal{H}_k$'s using the method of resurgence.
\end{rmk}

We end this introduction by summarizing the results in Theorems \ref{theorem: main-result-estimate} and \ref{theorem: quantizable-functions} as a comparison between classical and quantum objects. For a smooth function $f\in C^\infty(X)$, the associated quantum operator on the Hilbert spaces has an asymptotic expansion given by the higher Kostant--Souriau operators, which can be regarded as the quantum analogue of Taylor series expansion. (This is also consistent with the fact that a flat section under the Fedosov connection is a quantum Taylor expansion of a smooth function). For smooth functions $f$ which are symbols of quantizable functions, they are analogues of polynomials which have finite Taylor series expansions:
\begin{center}
	\renewcommand{\arraystretch}{1.5}
	\begin{tabular}{|c|c|}
		\hline
		{\color{red} Classical objects} & {\color{blue} Quantum objects}\\
		\hline
		smooth functions & operators on $\mathcal{H}_k$\\
		\hline
		polynomials & holomorphic differential operators on $\mathcal{H}_k$\\
		\hline
		terms in Taylor series expansion & higher Kostant--Souriau operators\\
		\hline
	\end{tabular}
\end{center}

\subsection*{Notations}
\begin{itemize}
	\item The same notation $\nabla$ denotes the naturally defined connections on various vector bundles such as the Weyl bundle, the pre-quantum line bundle $L$ and its tensor powers, etc, throughout this paper. 
	\item $D_{BT}$ denotes the (formal) Fedosov connection on the Weyl bundle whose associated star product is the Berezin--Toeplitz star product $\star_{BT}$. 
    $D_{BF,k}$ denotes the Fedosov connection on the level $k$ Bargmann--Fock sheaf $\mathcal{F}_k = \mathcal{W}_X\otimes_{\mathcal{O}_X}L^{\otimes k}$. 
	\item For every smooth function $f\in C^\infty(X)$, $P_{f,k,m}$ denotes the differential operator on $L^{\otimes k}$ corresponding to the weight $m$ component $(O_f)_m$ of the associated flat section. 
	\item $\W_{X,\C}:=\widehat{\Sym}TX_\C^*$ denotes the {\em complexified Weyl bundle} on $X$. 
	\item $C^\infty_{q,k}$ denotes the (sheaf of) level $k$ quantizable functions on $X$.
	\item We denote by $\lVert \cdot \rVert_{op}$ the operator norm of operators from $\mathcal{H}_k = H^0(X, L^{\otimes k})$ to $\Gamma(X, L^{\otimes k})$.
\end{itemize}

\subsection*{Acknowledgement} The first-named author would like to thank Martin Schlichenmaier for some early discussions via emails. The third-named author would like to thank Hang Xu and Zuoqin Wang for useful discussions on Berezin--Toeplitz operators. We are grateful to Leonid Polterovich for very useful comments and suggestions on an earlier draft of this paper. We also thank J\o rgen Ellegaard Andersen and Louis Ioos for their interest and encouragement.

K. Chan was supported by grants from the Hong Kong Research Grants Council (Project No. CUHK14305322, CUHK14305023 \& CUHK14302524). 
N. C. Leung was supported by grants from the Hong Kong Research
Grants Council (Project No. CUHK 14302224, CUHK 14305923 \& CUHK 14306322). 
Q. Li was supported by grants from National Natural Science Foundation of China (Project No. 12471061).

\section{Fedosov's approach to quantization}

In this section, we first give a quick review of the construction of a class of Fedosov connections on K\"ahler manifolds introduced in \cite{ChaLeuLi2022b}.  More details on these connections will be explained in Appendix \ref{appendix: Fedosov-connection-via-L-infinity-structure}. We will also introduce a polarized weight on the formal Weyl bundle. 

This type of Fedosov connections not only induce Wick type deformation quantization on K\"ahler manifolds, but also the Hilbert spaces of the geometric quantization on a pre-quantizable K\"ahler manifold $X$ via the Bargmann--Fock sheaf introduced in \cite{ChaLeuLi2021}. This is a module sheaf over the finite weight subbundle of $\W_{X,\C}[[\hbar]]$. There exists a Fedosov connection on this sheaf so that the space of flat sections is isomorphic to the Hilbert space $\mathcal{H}_k = H^0(X, L^{\otimes k})$. The Fedosov connections on the Weyl and the Bargmann-Fock bundles are compatible in a natural way.


\subsection{Fedosov deformation quantization on K\"ahler manifolds}
\

Recall that a {\em deformation quantization} of a symplectic manifold $(X, \omega)$ is a formal deformation of the commutative algebra $(C^{\infty }(X),\cdot)$ equipped with pointwise multiplication to a noncommutative one $( C^{\infty }( X)[[\hbar]] ,\star) $ equipped with a {\em star product} of the following form
$$
f\star g=fg+\sum_{i\geq 1}\hbar^i\cdot C_i(f,g),
$$
where each $C_i(-,-)$ is a bi-differential operator, so that the leading order of noncommutativity is the Poisson bracket $\{-,-\}$ associated to $\omega$, i.e.,
\begin{equation}\label{equation: Poisson-bracket}
	C_1(f,g)-C_1(g,f)= \left. \frac{d}{d\hbar}\left( f\star g - g\star f\right)\right\vert_{\hbar=0}
	= \{ f,g \}.
\end{equation}

A deformation quantization on a K\"ahler manifold $(X, J, \omega)$ is said to be of {\em Wick type} (also known as {\em separation of variables}) if all the bi-differential operators $C_l(f,g)$'s take only holomorphic derivatives of $f$ and only anti-holomorphic derivatives of $g$. It was shown by Karabegov in \cite{Kar1996} that to every Wick type star product, there is an associated closed formal $(1,1)$-form $-\frac{1}{\hbar} \omega +\alpha_0+ \alpha_1 \hbar + \alpha_2 \hbar^2 + \alpha_3 \hbar^3 + \cdots$, known as the {\em Karabegov form}, which gives rise to a one-to-one correspondence. 

In \cite{ChaLeuLi2022b}, it was shown that every Wick type deformation quantization on a K\"ahler manifold $X$ could be obtained by a Fedosov connection which arose from the quantization of the natural $L_\infty$ structure on $X$ \cite{Kap1999}. In particular, we can obtain the Berezin--Toeplitz star product $\star_{BT}$.
To recall the construction of Fedosov connections in \cite{ChaLeuLi2022b}, we first write the K\"ahler form $\omega$ on $X$ in local coordinates as
$$\omega=\omega_{i\bar{j}}dz^i\wedge d\bar{z}^j,$$
where we adopt the convention that $\omega^{\bar{k}i}\omega_{i\bar{j}}=\delta_{\bar{j}}^{\bar{k}}$ and $\omega^{i\bar{j}} = -\omega^{\bar{j}i}$.
We then consider the following {\em Weyl bundles} on $X$:
\begin{align*}\label{equation: Weyl-bundle}
	\W_{X}& := \widehat{\Sym}T^*X, \quad \overline{\W}_X:=\widehat{\Sym}\overline{T^*X}, \quad \W_{X,\C} :=\widehat{\Sym}T^*X_{\C}.
\end{align*}
To give explicit expressions of these bundles, we let $(z^1,\cdots, z^n)$ be a local holomorphic coordinate system on $X$, use $dz^i,d\bar{z}^j$'s to denote $1$-forms in $\A_X^\bullet$ and use $y^i,\bar{y}^j$ to denote sections of $T^*X_{\C}$, and also as generators of $\W_{X,\C}$. The K\"ahler form enables us to define a non-commutative fiberwise Wick product on $\W_{X,\C}[[\hbar]]$: 
\begin{equation}\label{equation: fiberwise-Wick-product}
	a\star b := \sum_{k\geq 0}\frac{\hbar^k}{k!}\cdot\omega^{i_1\bar{j}_1}\cdots\omega^{i_k\bar{j}_k}\cdot\frac{\partial^k a}{\partial y^{i_1}\cdots\partial y^{i_k}}\frac{\partial^k b}{\partial \bar{y}^{j_1}\cdots\partial \bar{y}^{j_k}}.
\end{equation}



The {\em symbol map}
\begin{equation}\label{equation: symbol-map}
	\sigma: \W_{X,\C}[[\hbar]]\rightarrow C^\infty(X)[[\hbar]].
\end{equation}
is defined by setting all $y^i,\bar{y}^j$'s to zero. 

\begin{defn}
	We will use the notation $\mathcal{W}_{p,q}$ to denote the component $\Sym^p T^*X\otimes_{\mathcal{C}^\infty_X}\Sym^q\overline{T^*X}$ of $\mathcal{W}_{X,\mathbb{C}}$; sections of this subbundle are said to be {\em of type $(p,q)$}. There are the following natural operators acting as derivations on $\A_X^\bullet(\mathcal{W}_{X,\mathbb{C}})$:
	\begin{align*}
		\delta^{1,0} a  = dz^i\wedge\frac{\partial a}{\partial y^i},\quad 
		\delta^{0,1}a  = d\bar{z}^j\wedge\frac{\partial a}{\partial\bar{y}^j},
	\end{align*}
	as well as
	\begin{align*}
		\delta^*a  = y^k\cdot \iota_{\partial_{z^k}}a+ \bar{y}^j\cdot \iota_{\partial_{\bar{z}^j}}a.
	\end{align*}
	We define the operator $\delta^{-1}$  by normalizing $\delta^{*}$ : 
	\begin{equation*}\label{equation: delta-1-0-inverse}
		\delta^{-1}:=\frac{\delta^*}{p_1+p_2+q_1+q_2}\ \text{on $\A_X^{p_1,q_1}(\mathcal{W}_{p_2,q_2})$},
	\end{equation*}
\end{defn}
 We also define the fiberwise de Rham differential as $\delta:=\delta^{1,0}+\delta^{0,1}$.

\begin{thm}[\cite{ChaLeuLi2022b}, Theorems 2.17 and 2.25]\label{theorem: Fedosov-connection}
	Let $\alpha = \sum_{i\geq 0}\hbar^i\alpha_i$ be a representative of a formal cohomology class in $H^2_{dR}(X)[[\hbar]]$ of type $(1,1)$. Then there exists $I_\alpha = I+ J_\alpha \in \A_X^{0,1}(\W_{X,\C}[[\hbar]])$, such that the connection 
	$$D_{\alpha} := \nabla-\delta+\frac{1}{\hbar}[I_\alpha, -]_{\star}$$
	is a Fedosov abelian connection. 
	The deformation quantization associated to the flat connection $D_{\alpha}$ is a Wick type star product with Karabegov form $K_{\alpha}$ given by $-\frac{1}{\hbar}\cdot\omega+\alpha$. 
\end{thm}
For the explicit expression of the term $I_\alpha$, we refer to \cite{ChaLeuLi2022b}. The connection $D_\alpha$ can also be written as 
$D_{\alpha}=\nabla+\frac{1}{\hbar}[\gamma_\alpha,-]_\star$.
The flatness of $D_{\alpha}$ is then equivalent to the following {\em Fedosov equation}:
\begin{equation}\label{equation: Fedosov-equation-gamma}
	R_\nabla+\nabla\gamma_\alpha+\frac{1}{\hbar}\gamma_\alpha\star\gamma_\alpha=-\omega+\hbar\cdot\alpha.
\end{equation}

In this paper,  we will focus on the Berezin--Toeplitz deformation quantization, whose Karabegov form is precisely given by
$$-\frac{1}{\hbar}\cdot\omega + \text{Ric}_X.$$ 
We will denote the associated Fedosov connection by $D_{BT}$ (where the subscript ``BT'' is short for ``Berezin--Toeplitz''), which induces a cochain complex
$$
\left(\A_X^\bullet(\W_{X,\C}[[\hbar]]), D_{BT}\right).
$$
We call this the {\em Fedosov complex}. 

The main result in \cite{Fed1994} states that there is a one-to-one correspondence between formal smooth functions on $X$ and flat sections of the Weyl bundle under Fedosov's abelian connection via the symbol map, which induces a star product on formal smooth functions.

\subsection{Polarized weight and flat sections associated to smooth functions}
\

In this subsection, we will first introduce a {\em polarized weight} on the Fedosov complex $\left(\A_X^\bullet(\W_{X,\C}[[\hbar]]), D_{BT}\right)$ as the sum of the polynomial degrees in $\hbar$ and the anti-holomorphic Weyl bundle $\overline{\W}_X$. For any smooth formal function $f\in C^\infty(X)[[\hbar]]$, we will explain how the flatness of $O_f$ implies equations connecting its weight components.

\subsubsection{The polarized weight on the Fedosov complex}
\

\begin{defn}\label{definition: polarized-weight}
	We define the {\em polarized weight} on the formal Weyl bundle $\W_{X,\C}[[\hbar]]$ by assigning weights on the generators:
	$$
	\lvert \hbar \rvert=1, \hspace{3mm} \lvert \bar{y}^j \rvert=1,\hspace{3mm}\lvert y^i \rvert=0.  
	$$
	This weight can be extended to $\A_X^\bullet(\W_{X,\C}[[\hbar]])$ by assigning weights $0$ to $\A_X^\bullet$. This weight is compatible with the fiberwise Wick product since $\lvert \hbar \rvert= \lvert y^i \rvert+ \lvert \bar{y}^j \rvert$.
\end{defn}

\begin{rmk}
	This weight is different from the one defined in \cite{Fed1994}. Throughout this paper, the notion {\em weight} on $\W_{X,\C}[[\hbar]]$ will mean the polarized weight in Definition \ref{definition: polarized-weight}. 
\end{rmk}

For any section $\alpha$ of the Weyl bundle $\W_{X,\C}[[\hbar]]$, let  $\alpha_i$ denote the weight $i$ component of $\alpha$.
We also introduce the following notation:
$$
\alpha_{\leq m}=\sum_{i=0}^m\alpha_i. 
$$

We now describe the weight decomposition of the Fedosov connection $D_{BT}$ on $\W_{X,\C}[[\hbar]]$. Recall that $D_{BT}$ is of the following explicit form (see Appendix \ref{subsection: Fedosov-quantization}, equation \eqref{equation: Fedosov-connection-Toeplitz-quantization}):
$$
D_{BT}=\nabla-\delta+\frac{1}{\hbar}[I_{BT},-]_\star=\nabla+\frac{1}{\hbar}[\gamma_{BT},-]_\star.
$$
A simple observation is that $D_{BT}$ can be decomposed into components of weight $-1$ and $0$; explicitly, the weight $-1$ component is given by
$$
-\delta^{0,1}
$$
(note that the differential forms $\A_X^\bullet$ have weight $0$), while the remaining terms give the weight $0$ component
$$
D_{BT}+\delta^{0,1}=\nabla-\delta^{1,0}+\frac{1}{\hbar}[I_{BT},-]_{\star}
$$
of $D_{BT}$.

The following picture illustrates the relation between the polarized weight on the Weyl bundle $\W_{X,\C}[[\hbar]]$ and the Fedosov connection $D_{BT}$. The horizontal coordinate denotes the degree of $\hbar$, and the vertical coordinate denotes the polynomial degree in $\overline{\W}_{X}$. Thus the weight of a monomial term in $\W_{X,\C}[[\hbar]]$ is the sum of these two coordinates.
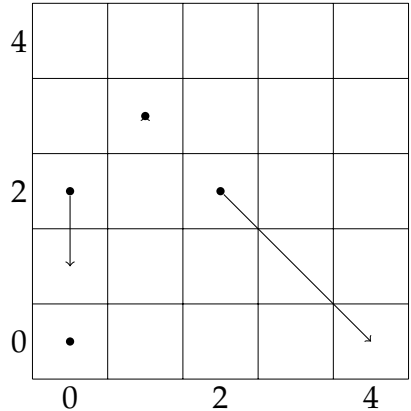
\begin{center}
	\begin{sseq}[entrysize=1cm]{0...4}{0...4} 
		\ssdropbull 
		\ssmoveto 0 2 \ssdropbull \ssvoidarrow{0}{-1} 
		\ssmoveto 2 2 \ssdropbull \ssvoidarrow{2}{-2}
		\ssmoveto 1 3 \ssdropbull \ssvoidarrow{0}{0}
	\end{sseq}
\end{center}
As to the Fedosov connection $D_{BT}$, first of all, its type $(1,0)$ part (with respect to the complex structure on $X$)
$$
D_{BT}^{1,0}=\nabla^{1,0}-\delta^{1,0}
$$
fixes both these two coordinates. Thus $D_{BT}^{1,0}$ fixes each position in the above diagram. Also, the term $-\delta^{0,1}$ in $D_{BT}$ is represented by the downward vertical arrows. All the other terms in $D_{BT}$ are of weight $0$, thus fixing the slope $-1$ line.

\subsubsection{Smooth functions and associated flat sections}
\

From Fedosov's original work \cite{Fed1994}, we know that there is a one-to-one correspondence between flat sections under a Fedosov connection and formal smooth functions
$$
\Gamma^{flat}(X, \W_{X,\C}[[\hbar]])\cong C^\infty(X)[[\hbar]]. 
$$
More explicitly, for every smooth function $f\in C^\infty(X)$, there exists a unique flat section $O_f$ of $\W_{X,\C}[[\hbar]]$ such that its symbol is exactly the function $f$:
$$
\sigma(O_f)=f.
$$


Since the $(1,0)$ component of $D_{BT}^{1,0}$ of the Fedosov connection preserves the polarized weight, we must have
\begin{equation}\label{equation: weight-component-annihilated-by-D-1-0}
	D_{BT}^{1,0}(O_f)_{m}=(\nabla^{1,0}-\delta^{1,0})(O_f)_m=0. 
\end{equation}
We have the following explicit relation between different weight components of $O_f$, which will play an essential role in the proofs of Lemma \ref{lemma: partial-bar-D-f-on-s} and Theorem \ref{theorem: main}.

\begin{lem}\label{lemma: iteration-relation-weight-components-flat-section}
	There is the following relation
	\begin{equation}\label{equation: iteration-relation-weight-components-flat-section}
		\begin{aligned}
			\delta^{0,1}\left((O_f)_{m+1}\right)=D_{BT}\left((O_f)_{\leq m}\right).
		\end{aligned}
	\end{equation}
\end{lem}
\begin{proof}
	Let $D_{BT}^{0,1}$ denote the $(0,1)$ component of the Fedosov connection, which is the sum of its weight $0$ and $-1$ components. In particular, its weight $-1$ component is exactly $-\delta^{0,1}$. Notice that the flatness of $O_f$ implis that 
	$$
	0=D_{BT}^{0,1}(O_f)=D_{BT}^{0,1}((O_f)_{\leq m})+D_{BT}^{0,1}((O_f)_{> m}),
	$$
	where $(O_f)_{> m} = \sum_{i=m+1}^\infty (O_f)_i$. It follows that 
	\begin{equation}\label{equation: iterative-relation-flat-section}
		D_{BT}((O_f)_{\leq m})=D_{BT}^{0,1}((O_f)_{\leq m})=\delta^{0,1}((O_f)_{m+1}).
	\end{equation}	
	Here in the first equality we have used the fact that each weight component of $O_f$ is annihilated by $D_{BT}^{1,0}=\nabla^{1,0}-\delta^{1,0}$ (equation \eqref{equation: weight-component-annihilated-by-D-1-0}). 
\end{proof}

\begin{rmk}
	Since each weight component of $O_f$ is a polynomial in $\hbar$, we can take the evaluation 
	$\hbar=\tfrac{\sqrt{-1}}{k}$ in equation \eqref{equation: iteration-relation-weight-components-flat-section}. This evaluation will play an essential role when we include geometric quantization in Fedosov's approach where $k$ will be the tensor power of the pre-quantum line bundle. 
\end{rmk}



\subsection{Bargmann--Fock sheaves and Hilbert spaces}
\

As a formal star product, the Wick product involves a formal variable $\hbar$. To couple it with geometric quantization, 
we take the {\em level $k$} evaluation $\hbar=\tfrac{\sqrt{-1}}{k}$ (the inverse of the tensor power of the pre-quantum line bundle, up to a scalar multiple) in the Wick product and define:
\begin{equation}\label{equation: level-k-Wick-product}
	\alpha\star_k\beta:=(\alpha\star\beta) \vert_{\hbar=\frac{\sqrt{-1}}{k}}.
\end{equation}
Similarly, we have the level $k$ Fedosov connection $D_{BT,k}$ (see equation \eqref{equation: level-k-Fedosov-connection} in Appendix \ref{appendix: Fedosov-connection-via-L-infinity-structure}).  
It is important to note that $\star_k$ and $D_{BT, k}$ are well-defined on the subbundle $$\widehat{\Sym}T^*X \otimes \Sym\overline{T^*X}[\hbar],$$
which we call the {\em finite weight subbundle} of $\W_{X,\C}[[\hbar]]$, as there will not be convergence issues related to power series in $\tfrac{1}{k}$ on this subbundle. 

The K\"ahler form on $X$ enables us to define the level $k$ fiberwise Bargmann--Fock action: a monomial in the finite weight subbundle acts as a differential operator on $\W_X=\widehat{\Sym}T^*X$  by the following {\em Wick ordering} formula:
\begin{equation}\label{equation: Wick-ordering-formula}
	\hbar^r\cdot y^{i_1}\cdots y^{i_m}\bar{y}^{j_1}\cdots\bar{y}^{j_l}\mapsto \left(\frac{\sqrt{-1}}{k}\right)^{l+r} \omega^{\bar{j}_1p_1}\cdots\omega^{\bar{j}_lp_l}\frac{\partial^l}{\partial y^{p_1} \cdots \partial y^{p_l}}\circ m_{y^{i_1}\cdots y^{i_m}}.
\end{equation}

This action denoted by $\circledast_k$ is compatible with $\star_k$ in the following sense: for any sections $\alpha_1,\alpha_2$ of the finite weight subbundle of $\W_{X,\C}[[\hbar]]$ and any section $s$ of $\W_X$,
\begin{equation}\label{equation: Bargmann-Fock-action-module}
	(\alpha_1 \star_k \alpha_2)\circledast_k s=\alpha_1\circledast_k(\alpha_2\circledast_k s).
\end{equation}
It makes $\W_X$ a sheaf of modules over the finite weight subbundle of $\W_{X,\C}[[\hbar]]$.

\begin{defn}[\cite{ChaLeuLi2023}, Definition 4.2]
	For every positive integer $k$,  we define the {\em level $k$ Bargmann--Fock sheaf} $\mathcal{F}_{L^{\otimes k}}$ by twisting $\W_X$ with the $k$-th tensor power of the prequantum line bundle $L$:
	$$
	\mathcal{F}_{L^{\otimes k}}:=\W_X\otimes_{\mathcal{O}_X}L^{\otimes k}.
	$$
\end{defn}
The Bargmann--Fock sheaf is a module over the finite weight subbundle of $\W_{X,\C}[[\hbar]]$ via the fiberwise Bargmann--Fock action. It was shown in \cite[Proposition 4.3]{ChaLeuLi2023} that $\mathcal{F}_{L^{\otimes k}}$ also admits a flat connection $D_{BF,k}$ given by
$$
D_{BF,k} = \nabla + \left(\frac{1}{\hbar}\cdot\gamma_{BT}\right)  \circledast_k-. 
$$
Here the subscript ``BF'' stands for ``Bargmann--Fock'', and $\gamma_{BT}$ is the same as in $D_{BT}$ in equation \eqref{equation: Fedosov-connection-Toeplitz-quantization}.  Here $\nabla$ denotes the tensor product of the Levi-Civita connection on $\W_{X}$ and the Chern connection on $L^{\otimes k}$. 


There is the following symbol map for the Bargmann--Fock sheaf:
\begin{equation}\label{equation: symbol-map-Bargmann-Fock}
	\sigma_{L^{\otimes k}}: \A_X^\bullet(\mathcal{F}_{L^{\otimes k}}) \rightarrow \A_X^\bullet(L^{\otimes k}),
\end{equation}
which is defined as the identity map for the differential form $\A_X^\bullet$ and $L^{\otimes k}$ part and by sending all $y^i$'s to zero.

Similar to the one-to-one correspondence between formal smooth functions and flat sections of the Weyl bundle, we have the following theorem.
\begin{thm}[\cite{ChaLeuLi2023}, Theorem 1.2]\label{theorem: flat-section-Hilbert-space}
	The symbol map $\sigma_{L^{\otimes k}}$ induces an isomorphism
	$$
	\Gamma^{flat}(X,\mathcal{F}_{L^{\otimes k}})\cong H^0(X, L^{\otimes k}) 
	$$
	from the space of flat sections of the Bargmann--Fock sheaf to the space of holomorphic sections of $L^{\otimes k}$.
\end{thm}
By applying equation \eqref{equation: Bargmann-Fock-action-module}, we can easily obtain the following compatibility of flat connections: for any section $\alpha$ of the finite weight subbundle of $\W_{X,\C}[[\hbar]]$ and any section $\gamma$ of $\mathcal{F}_{L^{\otimes k}}$,
\begin{equation}\label{equation: Fedosov-connection-compatibility}
	D_{BF,k}(\alpha\circledast_k \gamma)=D_{BT,k}(\alpha)\circledast_k\gamma+\alpha\circledast_k D_{BF,k}(\gamma). 
\end{equation}

For later computation, we will need an explicit expression for the flat section $O_s$ of $\mathcal{F}_{L^{\otimes k}}$ associated to a holomorphic section $s\in H^0(X, L^{\otimes k})$, as described in the following proposition.


\begin{prop}\label{proposition: expression-flat-section-Fock-bundle}
	The flat section $O_s$ associated to a holomorphic section $s\in H^0(X, L^{\otimes k})$ can be expressed as the sum
	\begin{equation}\label{equation: flat-section-Fock-bundle}
		O_s=\sum_{m\geq 0}(\tilde{\nabla}^{1,0})^m(s).
	\end{equation}
Here $\tilde{\nabla}^{1,0}:=\delta^{-1}\circ\nabla^{1,0}$.
\end{prop}
\begin{proof}
	 We consider $\tilde{O}_s:=\sum_{m\geq 0}(\tilde{\nabla}^{1,0})^m(s)$. Then it is obvious from the construction that 
	$$
	D_{BF,k}^{1,0}(\tilde{O}_s)=(\nabla^{1,0}-\delta^{1,0})(\tilde{O}_s)=0. 
	$$
	Thus both $O_s$ and $\tilde{O}_s$ are annihilated by $D_{BF,k}^{1,0}$, and they share the same symbol: $\sigma_{L^{\otimes k}}(O_s)=\sigma_{L^{\otimes k}}(\tilde{O}_s)=s$. There must be $O_s=\tilde{O}_s$, since otherwise we can take the term in $O_s-\tilde{O}_s$ of the lowest polynomial degree, which is not annihilated by $\delta^{1,0}$ and contradicts the fact that $D_{BF,k}^{1,0}(O_s-\tilde{O}_s)=0$.  
\end{proof}




\section{Higher Kostant-Souriau operators}

In this section, we define for each smooth function $f\in C^\infty(X)$ a series of differential operators acting on the Hilbert spaces $\mathcal{H}_k$ via the Bargmann--Fock action in the Fedosov quantization scheme. In Section \ref{subsection: Bargmann-Fock-higher-covariant-derivative}, we will identify the Bargmann-Fock action with certain differential-geometric differential operators known as {\em higher covariant derivatives}, of which we can establish norm estimates. 

Suppose the flat section associated to a function $f\in C^\infty(X)$ satisifes a finiteness condition, we call $O_f$ quantizable. 
In this case, the level $k$ Bargmann-Fock action $O_f\circledast_{k}-$ defines a holomorphic differential operator on the Hilbert space $\mathcal{H}_k$ via the identification $H^0(X, L^{\otimes k})\cong\Gamma^{\text{flat}}(X, \mathcal{F}_{L^{\otimes k}})$ and the compatibility between Fedosov connections on Weyl and Bargmann-Fock bundles. 

For a general smooth function $f\in C^\infty(X)$, we consider a truncation of $O_f$ using the polarized weight on $\W_{X,\C}$ in Definition \ref{definition: polarized-weight}, and define the higher Kostant-Souriau operators as the Bargmann-Fock action followed by the symbol map:
$$
P_{f,k,\leq m}:=\sigma_{L^{\otimes k}}\circ \left((O_f)_{\leq m}\circledast_k -\right).
$$

These Kostant-Souriau operators $P_{f,k,\leq m}$'s map holomorphic sections in $H^0(X, L^{\otimes k})$ to smooth sections $\Gamma(X, L^{\otimes k})$. By using the Fedosov equations for $O_f$ and $O_s$, we obtain an estimate of the norm of
\begin{equation}\label{equation: partial-bar-estimate}
	\bar{\partial}\left(P_{f,k,\leq m}(s)\right)
\end{equation}
in Proposition \ref{proposition: estimate-partial-bar-D-f-action-s} which shows that these Kostant-Souriau operators are "approximately" holomorphic.




\subsection{Bargmann--Fock action and higher covariant derivatives}\label{subsection: Bargmann-Fock-higher-covariant-derivative}
\

In this subsection, we show that given any monomial $\alpha$ in the Weyl bundle $\W_{X,\C}[[\hbar]]$ and any flat section $O_s$ of the Bargmann--Fock sheaf, the term $\sigma_{L^{\otimes k}}(\alpha\circledast_k O_s)$ can be expressed in terms of {\em higher covariant derivatives} of $s$ in differential geometry. We then give norm estimates of these operators when the domain is the space $\mathcal{H}_k$ of holomorphic sections. 

Here and in the rest of this paper, $\langle \cdot, \cdot \rangle$ denotes the naturally defined inner product of (smooth) sections of $L^{\otimes k}$, and $\lVert \cdot \rVert$ denotes the induced $L^2$-norm.

\subsubsection{Higher covariant derivatives}
\

We first give the definition of higher covariant derivatives for tensor fields.
\begin{lem-defn}\label{definition: higher-covariant-derivative}
	Let $E$ be a vector bundle on $X$ with a connection $\nabla^E$, and let $l \in \mathbb{N}$. Given any degree $l$ (not necessarily symmetric) tensor $G\in\Gamma(X, TX_{\C}^{\otimes l})$, we can define a {\em degree $l$ covariant derivative} on any section $s\in\Gamma(X, E)$ iteratively as follows. If $l = 0$ so that $G \in C^\infty(X)$, $\nabla_G^E(s) := G \cdot s$; if $l > 0$ and $G = X_1 \otimes G'$ for some $G' \in \Gamma(X, TX_{\C}^{\otimes (l-1)})$,
	$$
	\nabla_G^E(s):=\nabla_{X_1}^E\left(\nabla_{G'}^E(s)\right)-\nabla_{\nabla_{X_1}G'}^E(s).
	$$
	Moreover, $\nabla_G^E$ is a differential operator of order $l$ and is tensorial in $G$: for any $f \in C^\infty(X)$,
	$$
	\nabla_{f\cdot G}^E(s)=f\cdot\nabla_{G}^E(s).
	$$
\end{lem-defn}
\begin{proof}
	By induction on $l$, we will show that $\nabla_G$ is a well-defined differential operator of order $l$ and is tensorial in $G$. For $l=0$, these are the basic properties of multiplication by a function. For $l>0$, on one hand, considering the expression $(fX_1) \otimes G'$, there is by induction hypothesis
	$$
	\nabla_{f\cdot X_1}^E\left(\nabla_{G'}^E(s)\right)-\nabla_{\nabla_{f\cdot X_1}G'}^E(s)=f\cdot \left(\nabla_{X_1}^E\left(\nabla_{G'}^E(s)\right)-\nabla_{\nabla_{X_1}G'}^E(s)\right).
	$$
	On the other hand, considering $X_1 \otimes (fG')$, by the induction hypothesis there is
	\begin{align*}
		&\nabla_{X_1}^E\left(\nabla_{f\cdot G'}^E(s)\right)-\nabla_{\nabla_{X_1}(f\cdot G')}^E(s)\\
		=&\nabla_{X_1}^E\left(f\cdot\nabla_{ G'}^E(s)\right)-\nabla_{f\cdot\nabla_{X_1} G'}^E(s)-\nabla_{X_1(f)\cdot G'}^E(s)\\
		=&X_1(f)\cdot\nabla_{ G'}^E(s)+f\cdot\nabla_{X_1}^E\left(\nabla_{ G'}^E(s)\right)-\nabla_{f\cdot\nabla_{X_1} G'}^E(s)-X_1(f)\cdot\nabla_{G'}^E(s)\\
		=&f\cdot \left(\nabla_{X_1}^E\left(\nabla_{ G'}^E(s)\right)-\nabla_{\nabla_{X_1}G'}^E(s)\right).
	\end{align*}
	By the induction hypothesis again, the fact that $\nabla_G^E$ is a differential operator of order $l$ follows from an easy observation on the recursive formula.
\end{proof}
\begin{rmk}
	For $l=1$, we can easily check that the higher covariant derivative reduces to the standard covariant derivative.
\end{rmk}

\subsubsection{A norm estimate of higher covariant derivatives of holomorphic sections}
\quad\par
The goal of this subsection is to prove the following theorem about a norm estimate of higher covariant derivatives acting on holomorphic sections of $L^{\otimes k}$.

\begin{thm}\label{theorem: estimate-higher-covariant-derivative}
	Suppose $m \in \mathbb{N}$ and $G \in \Gamma(X, (TX)^{\otimes m})$. Then there exists $C = C_G > 0$ such that for all $k \in \mathbb{Z}^+$ and $s \in H^0(X, L^{\otimes k})$, we have
	\begin{equation*}
		\lVert \nabla_G s \rVert \leq C k^{\frac{m}{2}} \lVert s \rVert.
	\end{equation*}
\end{thm}

Our strategy in the proof of Theorem \ref{theorem: estimate-higher-covariant-derivative} is to turn the higher covariant derivatives into a sum of consecutive (standard) covariant derivatives, for which we prove a norm estimate by double induction. As a preparation, we state two a priori inequalities.

\begin{lem}
	Let $k \in \mathbb{Z}^+$, $s, s' \in \Gamma(X, L^{\otimes k})$ and $Y \in \Gamma(X, TX)$. Then
	\begin{equation}
		\label{equation: first-a-priori-inequality}
		\lvert \langle \nabla_Y s, s' \rangle \rvert \leq \lvert \langle \operatorname{div} (Y)  s, s' \rangle \rvert + \lvert \langle s, \overline{\partial}_{\overline{Y}} s' \rangle \rvert,
	\end{equation}
	where $\operatorname{div}(Y)$ is the divergence of $Y$ with respect to $\omega^n$, i.e. $\mathcal{L}_Y (\omega^n) = \operatorname{div}(Y) \omega^n$.
\end{lem}

On the other hand, for any $Y, Z \in \Gamma(X, TX)$, $k \in \mathbb{Z}^+$ and $s \in \Gamma(X, L^{\otimes k})$,
\begin{equation*}
	\overline{\partial}_{\overline{Y}} \nabla_Z s = \nabla_Z \overline{\partial}_{\overline{Y}} s + \nabla_{Z'} s - \overline{\partial}_{\overline{Y'}} s + \tfrac{k}{\sqrt{-1}} \omega(\overline{Y}, Z) s.
\end{equation*}
Here, $Y' = \overline{\partial}_{\overline{Z}} Y$ and $Z' = \overline{\partial}_{\overline{Y}} Z$, whence $[\overline{Y}, Z] = Z' - \overline{Y'}$. 
It implies the following lemma.

\begin{lem}
	Let $l \in \mathbb{N}$, $k \in \mathbb{Z}^+$, $s, s' \in \Gamma(X, L^{\otimes k})$ and $Z_1, ..., Z_{l+1} \in \Gamma(X, TX)$. Then for all $i \in \mathbb{N}$ with $1 \leq i \leq l$, we have
	\begin{equation}
		\label{equation: second-a-priori-inequality}
		\begin{aligned}
			& \lvert \langle s, \nabla_{Z_1} \cdots \nabla_{Z_{l-i-1}} \overline{\partial}_{\overline{Z_{l-i}}} \nabla_{Z_{l-i+1}} \cdots \nabla_{Z_l} s' \rangle \rvert\\
			\leq & \lvert \langle s, \nabla_{Z_1} \cdots \nabla_{Z_{l-i-1}} \nabla_{Z_{l-i+1}} \overline{\partial}_{\overline{Z_{l-i}}} \nabla_{Z_{l-i+2}} \cdots \nabla_{Z_l} s' \rangle \rvert\\
			& + \lvert \langle s, \nabla_{Z_1} \cdots \nabla_{Z_{l-i-1}} \nabla_{Z'_{l-i+1}} \nabla_{Z_{l-i+2}} \cdots \nabla_{Z_l} s' \rangle \rvert\\
			& + \lvert \langle s, \nabla_{Z_1} \cdots \nabla_{Z_{l-i-1}} \overline{\partial}_{\overline{Z'_{l-i}}} \nabla_{Z_{l-i+2}} \cdots \nabla_{Z_l} s' \rangle \rvert\\
			& + \lvert \langle s, \nabla_{Z_1} \cdots \nabla_{Z_{l-i-1}} \tfrac{k}{\sqrt{-1}} \omega(\overline{Z_{l-i}}, Z_{l-i+1}) \nabla_{Z_{l-i+2}} \cdots \nabla_{Z_l} s' \rangle \rvert,
		\end{aligned}
	\end{equation}
	where $Z'_{l-i+1} = \overline{\partial}_{\overline{Z_{l-i}}} Z_{l-i+1}$ and $Z'_{l-i} = \overline{\partial}_{\overline{Z_{l-i+1}}} Z_{l-i}$.
\end{lem}

Now, we state a more general statement to be proved.

\begin{thm}
	\label{theorem: estimate-consecutive-covariant-derivative}
	Let $m \in \mathbb{N}$.
	\begin{enumerate}
		\item For all $l \in \mathbb{N}$ with $l \leq m$, the following statement, denoted by $P_m(l)$, holds:\\
		For all $Y_1, ..., Y_m, Z_1, ..., Z_l \in \Gamma(X, TX)$, there exists $C = C_{Y_1, ..., Y_m, Z_1, ..., Z_l} > 0$ such that for all $k \in \mathbb{Z}^+$ and $s \in H^0(X, L^{\otimes k})$,
		\begin{equation*}
			\lvert \langle \nabla_{Y_1} \cdots \nabla_{Y_m} s, \nabla_{Z_1} \cdots \nabla_{Z_l} s \rangle \rvert \leq C \cdot k^m \cdot \lVert s \rVert^2.
		\end{equation*}
		\item For all $i, l \in \mathbb{N}$ with $i \leq l \leq m$, the following statement, denoted by $P_m(l, i)$, holds:\\
		For all $Y_1, ..., Y_m, Z_0, ..., Z_l \in \Gamma(X, TX)$, there exists $C = C_{Y_1, ..., Y_m, Z_0, ..., Z_l} > 0$ such that for all $k \in \mathbb{Z}^+$ and $s \in H^0(X, L^{\otimes k})$,
		\begin{equation*}
			\lvert \langle \nabla_{Y_1} \cdots \nabla_{Y_m} s, \nabla_{Z_0} \cdots \nabla_{Z_{l-i-1}} \overline{\partial}_{\overline{Z_{l-i}}} \nabla_{Z_{l-i+1}} \cdots \nabla_{Z_l} s \rangle \rvert \leq C \cdot k^{m+1} \cdot \lVert s \rVert^2.
		\end{equation*}
	\end{enumerate}
\end{thm}

It is evident that $P_0(0)$ is true. Furthermore, for all $m, l \in \mathbb{N}$ with $l \leq m$, the statement $P_m(l, 0)$ holds. This is because, by the holomorphicity of $s$, we have
\begin{equation*}
	\lvert \langle \nabla_{Y_1} \cdots \nabla_{Y_m} s, \nabla_{Z_0} \cdots \nabla_{Z_{l-1}} \overline{\partial}_{\overline{Z_l}} s \rangle \rvert = 0.
\end{equation*}
Therefore, Theorem \ref{theorem: estimate-consecutive-covariant-derivative} holds when $m = 0$. From now on, assume that Theorem \ref{theorem: estimate-consecutive-covariant-derivative} holds for some $m \in \mathbb{N}$. We then continue our proof by induction, dividing the rest of our proof into two lemmas.

\begin{lem}
	\label{lemma: estimate-consecutive-derivative}
	For all $l \in \mathbb{N}$ with $l \leq m + 1$, $P_{m+1}(l)$ holds.
\end{lem}
\begin{proof}
	Fix $Y_1, ..., Y_{m+1}, Z_1, ..., Z_l \in \Gamma(X, TX)$, $k \in \mathbb{Z}^+$ and $s \in H^0(X, L^{\otimes k})$.
	\begin{enumerate}
		\item Suppose $l \leq m$. Then by the induction hypothesis that $P_m(l)$ and $P_m(l, l)$ hold,
		\begin{align*}
			\lvert \langle \operatorname{div} (Y_1)  \nabla_{Y_2} \cdots \nabla_{Y_{m+1}} s, \nabla_{Z_1} \cdots \nabla_{Z_l} s \rangle \rvert = & O(k^m) \lVert s \rVert^2,\\
			\lvert \langle \nabla_{Y_2} \cdots \nabla_{Y_{m+1}} s, \overline{\partial}_{\overline{Y_1}} \nabla_{Z_1} \cdots \nabla_{Z_l} s \rangle \rvert = & O(k^{m+1}) \lVert s \rVert^2.
		\end{align*}
		Therefore, by equation (\ref{equation: first-a-priori-inequality}), $P_{m+1}(l)$ holds.
		\item Suppose $l = m + 1$. Since we have already established that $P_{m+1}(m)$ holds in the previous case, it follows that
		\begin{equation*}
			\lvert \langle \operatorname{div} (Y_1)  \nabla_{Y_2} \cdots \nabla_{Y_{m+1}} s, \nabla_{Z_1} \cdots \nabla_{Z_{m+1}} s \rangle \rvert = O(k^{m+1}) \lVert s \rVert^2.
		\end{equation*}
		Note that $\lvert \langle \nabla_{Y_2} \cdots \nabla_{Y_{m+1}} s, \nabla_{Z_1} \cdots \nabla_{Z_{m+1}} \overline{\partial}_{\overline{Y_1}} s \rangle \rvert = 0$ and for all $1 \leq i \leq m+1$,
		\begin{align*}
			\lvert \langle \nabla_{Y_2} \cdots \nabla_{Y_{m+1}} s, \nabla_{Z_1} \cdots \nabla_{Z_{i-1}} \nabla_{Z'_i} \nabla_{Z_{i+1}} \cdots \nabla_{Z_{m+1}} s \rangle \rvert = & O(k^{m+1}) \lVert s \rVert^2,\\
			\lvert \langle \nabla_{Y_2} \cdots \nabla_{Y_{m+1}} s, \nabla_{Z_1} \cdots \nabla_{Z_{i-1}} \overline{\partial}_{\overline{Y'_i}} \nabla_{Z_{i+1}} \cdots \nabla_{Z_{m+1}} s \rangle \rvert = & O(k^{m+1}) \lVert s \rVert^2,\\
			\lvert \langle \nabla_{Y_2} \cdots \nabla_{Y_{m+1}} s, \nabla_{Z_1} \cdots \nabla_{Z_{i-1}} \tfrac{k}{\sqrt{-1}} \omega(\overline{Y_1}, Z_i) \nabla_{Z_{i+1}} \cdots \nabla_{Z_{m+1}} s \rangle \rvert = & O(k^{m+1}) \lVert s \rVert^2,
		\end{align*}
		because $P_{m+1}(m)$, $P_m(m, m+1-i)$ and $P_m(m)$ hold. Here, $Z'_i = \overline{\partial}_{\overline{Y_1}} Z_i$ and $Y'_i = \overline{\partial}_{\overline{Z_i}} Y_1$. Then by equation \eqref{equation: second-a-priori-inequality} and induction,
		\begin{equation*}
			\lvert \langle \nabla_{Y_2} \cdots \nabla_{Y_{m+1}} s, \overline{\partial}_{\overline{Y_1}} \nabla_{Z_1} \cdots \nabla_{Z_{m+1}} s \rangle \rvert = O(k^{m+1}) \lVert s \rVert^2.
		\end{equation*}
		Eventually, by equation (\ref{equation: first-a-priori-inequality}), $P_{m+1}(m+1)$ holds.
	\end{enumerate}
\end{proof}

\begin{lem}
	For all $i, l \in \mathbb{N}$ with $i \leq l \leq m + 1$, $P_{m+1}(l, i)$ holds.
\end{lem}
\begin{proof}
	We will argue by induction on $i$, knowing that $P_{m+1}(l, 0)$ holds for all $l \in \mathbb{N}$ with $l \leq m+1$. Now, assume $i \in \mathbb{N}$ with $i \leq m$ such that for all $l \in \mathbb{N}$ with $i \leq l \leq m+1$, $P_{m+1}(l, i)$ holds.\par
	Consider any $l \in \mathbb{N}$ with $i + 1 \leq l \leq m+1$, $Y_1, ..., Y_{m+1}, Z_0, ..., Z_l \in \Gamma(X, TX)$, $k \in \mathbb{Z}^+$ and $s \in H^0(X, L^{\otimes k})$. Let $Z'_{l-i} = \overline{\partial}_{\overline{Z_{l-i-1}}} Z_{l-i}$ and $Z'_{l-i-1} = \overline{\partial}_{\overline{Z_{l-i}}} Z_{l-i-1}$. By Lemma \ref{lemma: estimate-consecutive-derivative}, $P_{m+1}(l)$ and $P_{m+1}(l-1)$ hold. Thus,
	\begin{align*}
		\lvert \langle \nabla_{Y_1} \cdots \nabla_{Y_{m+1}} s, \nabla_{Z_0} \cdots \nabla_{Z_{l-i-2}} \nabla_{Z'_{l-i}} \nabla_{Z_{l-i+1}} \cdots \nabla_{Z_l} s \rangle \rvert \leq & O(k^{m+1}) \lVert s \rVert^2,\\
		\lvert \langle \nabla_{Y_1} \cdots \nabla_{Y_{m+1}} s, \nabla_{Z_0} \cdots \nabla_{Z_{l-i-2}} \tfrac{k}{\sqrt{-1}} \omega(\overline{Z_{l-i-1}}, Z_{l-i}) \nabla_{Z_{l-i+1}} \cdots \nabla_{Z_l} s \rangle \rvert = & O(k^{m+2}) \lVert s \rVert^2.
	\end{align*}
	By the induction hypothesis that $P_{m+1}(l, i)$ and $P_{m+1}(l-1, i)$ hold,
	\begin{align*}
		\lvert \langle \nabla_{Y_1} \cdots \nabla_{Y_{m+1}} s, \nabla_{Z_0} \cdots \nabla_{Z_{l-i-2}} \nabla_{Z_{l-i}} \overline{\partial}_{\overline{Z_{l-i-1}}} \nabla_{Z_{l-i+1}} \cdots \nabla_{Z_l} s \rangle \rvert = & O(k^{m+2}) \lVert s \rVert^2,\\
		\lvert \langle \nabla_{Y_1} \cdots \nabla_{Y_{m+1}} s, \nabla_{Z_0} \cdots \nabla_{Z_{l-i-2}} \overline{\partial}_{\overline{Z'_{l-i-1}}} \nabla_{Z_{l-i+1}} \cdots \nabla_{Z_l} s \rangle \rvert = & O(k^{m+2}) \lVert s \rVert^2.
	\end{align*}
	Finally, by (\ref{equation: second-a-priori-inequality}), $\lvert \langle \nabla_{Y_1} \cdots \nabla_{Y_{m+1}} s, \nabla_{Z_0} \cdots \nabla_{Z_{l-i-2}} \overline{\partial}_{\overline{Z_{l-i-1}}} \nabla_{Z_{l-i}} \cdots \nabla_{Z_l} s \rangle \rvert = O(k^{m+2}) \lVert s \rVert^2$.
\end{proof}

\begin{lem}\label{lemma: decomposition-of-tensor-G}
	For any $m\geq 1$, and any tensor $G\in\Gamma(X, TX^{\otimes m})$ can be written as a finite sum of tensors of the following form:
	$$
	X_1\otimes \cdots\otimes X_m,
	$$
	where  each $X_i\in\Gamma(X, TX)$. 
\end{lem}
\begin{proof}
	This is obviously true locally, and also holds globally on $X$ by applying a partition of unity. 
\end{proof}

By a simple induction on the degree of the tensor $G$, we can prove the following:
\begin{lem}\label{lemma: higher-covariant-derivative-consecutive-derivative}
	Let $m>0$ and $G\in\Gamma(X, TX^{\otimes m})$. The higher covariant derivative $\nabla_G$ can be written as a finite sum of consecutive covariant derivatives of the following form:  
	$$
	\nabla_{X_1}\circ\cdots\circ \nabla_{X_l},
	$$
	where $1 \leq l \leq m$ and all these $X_i$'s are global tangent vector fields on $X$ of type $(1,0)$.
\end{lem}
\begin{proof}
	We do induction on $m$. For $m=1$, this is obviously true. For $m\geq 2$, by Lemma \ref{lemma: decomposition-of-tensor-G}, we can assume without loss of generality that $G=Y_1\otimes\cdots\otimes Y_m$.  By  Definition \ref{definition: higher-covariant-derivative}, there is:
	$$
	\nabla_G(s)=\nabla_{Y_1}\circ\nabla_{Y_2\otimes\cdots\otimes Y_m}(s)-\nabla_{\nabla_{Y_1}(Y_2\otimes\cdots\otimes Y_m)}(s).
	$$ 
	and the statement holds for $m$ by the induction hypothesis. 
\end{proof}

\begin{proof}[Proof of Theorem \ref{theorem: estimate-higher-covariant-derivative}]
	When $m = 0$, as $G$ is a smooth function on a compact manifold $X$, its norm is bounded above by a constant $C > 0$, whence $\lVert \nabla_G s \rVert \leq C \lVert s \rVert$. When $m > 0$, by Lemma \ref{lemma: higher-covariant-derivative-consecutive-derivative}, $\nabla_G s$ is a finite sum of consecutive covariant derivatives of the form $\nabla_{X_1} \cdots \nabla_{X_l} s$, where $1 \leq l \leq m$ and $X_1, ..., X_l \in \Gamma(X, TX)$. In particular, this identification is independent of $k$, and hence by Theorem \ref{theorem: estimate-consecutive-covariant-derivative}, there exists $C > 0$ which is independent of $k$ and $s$ such that $\lVert \nabla_G s \rVert \leq C k^{\frac{m}{2}} \lVert s \rVert$.
\end{proof}

\subsubsection{Bargmann--Fock action and higher covariant derivatives}
\

Recall that the isomorphism in Theorem \ref{theorem: flat-section-Hilbert-space} implies that for any holomorphic section $s\in\mathcal{H}_k=H^0(X, L^{\otimes k})$, there is a canonically associated flat section of the level $k$ Bargmann--Fock sheaf, which we denote by $O_s$.

\begin{defn}\label{definition: differential-operators-finite-weight-section-Weyl}
	For any section $\alpha$ of the finite weight subbundle of $\W_{X,\C}[[\hbar]]$, we define the associated differential operator
	\begin{equation*}
		P_{\alpha,k}: \mathcal{H}_k=H^0(X, L^{\otimes k})\rightarrow \Gamma(X, L^{\otimes k})
	\end{equation*}
	as the composition of Bargmann--Fock action and the symbol map for all quantum levels $k>0$.  Explicitly, 
	\begin{equation}\label{equation: differential-operator-Bargmann-Fock}
		P_{\alpha,k}(s):=\sigma_{L^{\otimes k}} (\alpha\circledast_k O_s).
	\end{equation}
\end{defn}
\begin{rmk}
	It is not difficult to see from the definition of the Bargmann--Fock action in equation \eqref{equation: Wick-ordering-formula} that if $\alpha\in\Gamma(X, \hbar^r \cdot \Sym^l T^*X\otimes\Sym^m\overline{T^*X})$, where $l>m$, then the associated differential operator $P_{\alpha,k}$ must vanish for all $k>0$. 
\end{rmk}


The main result of this subsection is the following identification of the Bargmann--Fock action with higher covariant derivatives on $O_s$.
\begin{prop}\label{proposition: Bargmann-Fock-action-higher-covariant-derivative}
	Let $\alpha$ be a section of the Weyl bundle which lives in $\Sym^l T^*X\otimes\Sym^m\overline{T^*X}$ where $m\geq l$. Then there exists a tensor field $G\in\Gamma(X, TX^{\otimes (m-l)})$ depending only on $\alpha$ but independent of the level $k$, such that for all $k>0$, we have
	$$
	P_{\alpha,k}= \left( \frac{\sqrt{-1}}{k} \right)^m \cdot\nabla_{G}: \mathcal{H}_k\rightarrow\Gamma(X, L^{\otimes k}).  
	$$
	In other words, the differential operators defined via Bargmann--Fock action can be expressed as a {\em higher covariant derivative}. 
\end{prop}
\begin{proof}
	The statement is trivially true for $m = 0$. From now on, assume $m > 0$. Let us write $\alpha$ locally as
	\begin{equation}\label{equation: monomial-Weyl-bundle}
		\alpha=\alpha_{i_1\cdots i_l, \bar{j}_1\cdots \bar{j}_m}y^{i_1}\cdots y^{i_l}\bar{y}^{j_1}\cdots\bar{y}^{j_m}.
	\end{equation}

	Fix $s \in \mathcal{H}_k$. We will need the following formula for $O_s$ as in Proposition \ref{proposition: expression-flat-section-Fock-bundle}:
	$$
	O_s=\sum_{p\geq 0}\left(\tilde{\nabla}^{1,0}\right)^p(s). 
	$$
	Let $(O_s)_m:=\left(\tilde{\nabla}^{1,0}\right)^m(s)$ denote the term in $O_s$ of polynomial degree $m$. We first prove the case where $l=0$. In this case, the tensor field $G$ is of the following form:
	$$
	G= \alpha_{\bar{j_1}\cdots\bar{j}_m}\omega^{\bar{j}_1i_1}\cdots\omega^{\bar{j}_mi_m}\partial_{z^{i_1}}\cdots\partial_{z^{i_m}}=\frac{1}{m!}\alpha_{\bar{j_1}\cdots\bar{j}_m}\omega^{\bar{j}_1i_1}\cdots\omega^{\bar{j}_mi_m}\partial_{z^{i_1}}\otimes\cdots\otimes\partial_{z^{i_m}}.
	$$
	In the second equality, we have used the symmetry of the coefficients with respect to the indices.  Here the term $\left( \tfrac{\sqrt{-1}}{k} \right)^m$ arises from the assignment $\bar{y}^{j_p}$ to $\frac{\sqrt{-1}}{k}\cdot\omega^{\bar{j_p}i_p}\cdot\partial_{z^{i_p}}$.
	
	Without loss of generality, we can assume that $G_m=X_1\otimes\cdots\otimes X_m$, for vector fields $X_1,\cdots X_m$ on $X$ of type $(1,0)$. We claim that it is enough to prove the following equality:
	$$
	\nabla_{G_m}(s)=\langle G_m, (O_s)_m\rangle.
	$$
	Here $\langle -,-\rangle$ denotes the natural pairing between tensors. 
	We do induction on $m$: for $m=1$ this follows from the definition of covariant derivatives and the explicit formula of $O_s$. Now, consider $G_{m+1} = X_1 \otimes \cdots \otimes X_{m+1}$ and let $\tilde{G}_m = X_2 \otimes \cdots \otimes X_{m+1}$. By the inductive definition of {\em higher covariant derivatives}, there is
	\begin{align*}
		\nabla_{G_{m+1}}(s) :=&\nabla_{X_1}\left(\nabla_{\tilde{G}_m}(s)\right)-\nabla_{\nabla_{X_1} \tilde{G}_m}(s)\\
		=&\nabla_{X_1}\left(\langle \tilde{G}_m, (O_s)_m\rangle\right)-\langle \nabla_{X_1}\tilde{G}_m, (O_s)_m\rangle\\
		=&\langle \tilde{G}_m, \nabla_{X_1}(O_s)_m\rangle\\
		=&\langle X_1\otimes \tilde{G}_m, (\delta^{1,0})^{-1}\circ\nabla^{1,0}(O_s)_m\rangle\\
		=&\langle G_{m+1}, (O_s)_{m+1}\rangle.
	\end{align*}
	Here in the second step we have used the inductive hypothesis for $\tilde{G}_m$ and $\nabla_{X_1}\tilde{G}_m$, and in the last step we have used the following equality:
	$$
	(O_s)_{m+1}=(\delta^{1,0})^{-1}\circ\nabla^{1,0}(O_s)_m. 
	$$
	and also the symmetry property described in Lemma \ref{lemma: symmetry-property-flat-section}. 
	
	For the general case, we can assume that $m\geq l$ in $\alpha$ as otherwise the output will vanish after taking the symbol.  Notice that in the definition of Bargmann--Fock action, we will first multiply those $y^i$'s in $\alpha$ to $O_s$ and then take derivatives associated to $\bar{y}^j$'s.
	
	We claim that after taking the symbol, the Bargmann--Fock action can be written as a pure differentiation. More precisely, there exists $\beta=\beta_{j_1\cdots j_{m-l}}\bar{y}^{j_1}\cdots\bar{y}^{j_{m-l}}$, such that 
	$$
	\sigma_{L^{\otimes k}}\left(\alpha\circledast_k O_s\right)=\sigma_{L^{\otimes k}}\left(\beta\circledast_k O_s\right).
	$$
	Thus it is enough to prove for the cases where $l=0$. 
\end{proof}

\begin{lem}\label{lemma: symmetry-property-flat-section}
	For any holomorphic section $s\in H^0(X, L^{\otimes k})$, and any vector fields $X_1,\cdots, X_{l+1}$ on $X$ of type $(1,0)$, we have
	$$
	\langle X_2\otimes\cdots\otimes X_{l+1},\nabla_{X_1}\circ(\tilde{\nabla}^{1,0})^l(s)\rangle=\langle X_1\otimes\cdots \otimes X_{l+1}, (\tilde{\nabla}^{1,0})^{l+1}s\rangle. 
	$$
\end{lem}
\begin{proof}

	Let $\beta_l:=\nabla^{1,0}\circ(\tilde{\nabla}^{1,0})^l(s)=\beta_{i,j_1\cdots j_l}dz^i\otimes y^{j_1}\cdots y^{j_l}$, we claim that the index $i$  in the subscript of $\beta$ is symmetric with the other indices $j_1,\cdots, j_l$. (Here each $\beta_{i,j_1\cdots j_l}$ is a local smooth section of $L^{\otimes k}$.)
	
	To prove this claim, a simple observation is that this desired symmetry is equivalent to the vanishing:
	$$
	\delta^{1,0}\left(\beta_l\right)=0. 
	$$
	We will also need to use the following commutator relation:
	$$
	[\delta^{1,0},\nabla]=0,
	$$
	which follows from the symmetric property of the Christoffel symbol $\Gamma_{ij}^k=\Gamma_{ji}^k$.

	We will do induction on $l$. For $l=1$, this follows from the fact that $(\nabla^{1,0})^2=0$. In general there is
	\begin{align*}
		\delta^{1,0}\left(\nabla^{1,0}\circ(\tilde{\nabla}^{1,0})^{l+1}(s)\right)=&	\nabla^{1,0}\left(\delta^{1,0}\circ(\tilde{\nabla}^{1,0})^{l+1}(s)\right)\\
		=&\nabla^{1,0}\left(\delta^{1,0}\circ(\delta^{1,0})^{-1}\circ\nabla^{1,0}\circ(\tilde{\nabla}^{1,0})^l(s)\right)\\
		=&\nabla^{1,0}\left((id-(\delta^{1,0})^{-1}\circ\delta^{1,0})\circ\nabla^{1,0}\circ(\tilde{\nabla}^{1,0})^l(s)\right)\\
		=&\nabla^{1,0}\left(\beta_l-(\delta^{1,0})^{-1}\circ\delta^{1,0}(\beta_l)\right)\\
		=&\nabla^{1,0}(\beta_l)\\
		=&0.
	\end{align*}
	Here in the third equality we have used the following equation to obtain that $\delta^{1,0}\circ(\delta^{1,0})^{-1}(\beta_l)=\beta_l-(\delta^{1,0})^{-1}\circ\delta^{1,0}(\beta_l)$:
	$$	\delta^{1,0}\circ(\delta^{1,0})^{-1}+(\delta^{1,0})^{-1}\circ\delta^{1,0}=\text{id}-\pi^{1,0}.$$
	Here $\pi^{1,0}:\A_X^\bullet\otimes\mathcal{F}_k\rightarrow\A_X^{0,\bullet}\otimes L^{\otimes k}$ denotes  
	the map setting $y^i, dz^i$'s to $0$ and identity for $d\bar{z}^j$'s and sections of $L^{\otimes k}$, and it is clear that $\pi^{1,0}(\beta_l)=0$ for all $l\geq 1$. In the fifth equality we have used the induction hypothesis $\delta^{1,0}(\beta_l)=0$. And the last equality follows from the fact that $(\nabla^{1,0})^2=0$. 
	
	It follows from this claim that 
	$$
	(\tilde{\nabla}^{1,0})^{l+1}(s)=\beta_{i,j_1\cdots j_l}y^i\otimes y^{j_1}\otimes\cdots\otimes y^{j_l}.
	$$
	And there is
	\begin{align*}
		\langle X_2\otimes\cdots\otimes X_{l+1},\nabla_{X_1}\circ(\tilde{\nabla}^{1,0})^l(s)\rangle	=&\beta_{i,j_1\cdots j_l}\cdot X_1(dz^i)\cdot \langle X_2\otimes\cdots\otimes X_{l+1},y^{j_1}\otimes\cdots\otimes y^{j_l}\rangle\\
		=&\beta_{i,j_1\cdots j_l}\cdot\langle X_1\otimes X_2\otimes\cdots\otimes X_{l+1}, y^i\otimes y^{j_1}\otimes\cdots\otimes y^{j_l}\rangle.
	\end{align*}
\end{proof}

\subsection{Quantizable functions and holomorphic differential operators}\label{subsection: quantizable-function}
\

In this subsection, we first review the notion of quantizable functions as introduced in \cite{ChaLeuLi2023}. There are two types of quantizable functions: formal and non-formal quantizable functions, which are closely related. We begin with the formal one.
\begin{defn}[\cite{ChaLeuLi2023}, Definition 2.13]\label{definition: formal-quantizable-function}
	A formal function $f\in C^\infty(X)[[\hbar]]$ is called a {\em formal quantizable function} if its associated flat section $O_f$ lives in the finite weight subbundle of $\W_{X,\C}[[\hbar]]$. A formal quantizable function $f$ is said to be of {\em order $m$} if the highest weight component of $O_f$ is of weight $m$.
\end{defn}

\begin{eg}\label{example: holomorphic-function-quantizable}
	Formal quantizable functions on $\mathbb{C}^n$ are explicitly given by
	$$
	\mathcal{O}_{\C^n}[\bar{z}^1,\cdots,\bar{z}^n,\hbar],
	$$
	i.e., polynomials in $\bar{z}^1,\cdots,\bar{z}^n$ and $\hbar$ with coefficients holomorphic functions on $\C^n$. 
\end{eg}
\begin{eg}\label{example: holomorphic-functions}
	On a K\"ahler manifold $X$, every (local) holomorphic function $f$ is a formal quantizable function for {\em any} Karabegov form. Explicitly, the flat section associated to $f$ is given by 
	$$
	O_f=\sum_{k\geq 0}(\tilde{\nabla}^{1,0})^k(f),
	$$
	thus  $O_f\in\left(\W_{X,C}[[\hbar]]\right)_0=\W_X$, making $f$ an order $0$ formal quantizable function.
\end{eg}
\begin{rmk}\label{remark: quantizable functions}
	In the picture of the Kostant--Souriau geometric quantization, there was also a notion of quantizable functions, defined as functions whose Hamiltonian vector fields are polarization-preserving. On $\C^n$, it is easy to check that these functions are at most of degree $1$ in $\bar{z}^1,\cdots,\bar{z}^n$. Thus quantizable functions defined here and in \cite{ChaLeuLi2023} are much more general than those in the Kostant--Souriau geometric quantization. 
\end{rmk}

\begin{eg}
	In \cite{LeuLiMa2024}, the (global) quantizable functions of order one on a K\"ahler manifold $X$ were carefully studied. 
	It was shown in \cite[Proposition 3.1]{LeuLiMa2024} that given a smooth function $f_0\in C^\infty(X)$ such that 
	\begin{equation}\label{equation: f-0}
		\nabla^{0,1}\left(\frac{\partial f_0}{\partial\bar{z}^j}\bar{y}^j\right)=0.
	\end{equation}
	Then $f_0-\hbar \cdot \frac{1}{\sqrt{-1}}\Delta f_0$ is a first order formal quantizable function. 
	In particular, the Hamiltonian function of a vector field $V$ preserving the complex polarization satisfies equation \eqref{equation: f-0}. In this case, the function $f_0-\hbar \cdot \frac{1}{\sqrt{-1}}\Delta f_0$ is called a {\em quantum moment map}. This shows that quantum moment maps give examples of first order quantizable functions. 
\end{eg}

For every fixed level $k \geq 1$, there are non-formal Fedosov connections $D_{BT,k}$ (see Appendix \ref{subsection: non-formal-Fedosov-connection}), with which we can also define the notion of non-formal quantizable functions; see \cite{ChaLeuLi2023} for more details.
\begin{defn}[\cite{ChaLeuLi2023}, Definition 2.20] \label{definition: quantizable functions}
	A flat section $\alpha\in\Gamma(X, \W_X\otimes\Sym\overline{T^*X})$ under the non-formal Fedosov connection $D_{BT,k}$ is called a {\em (non-formal) quantizable function} of level $k$.
	We denote the space of these functions by $C_{q,k}^\infty(X)$. 
\end{defn}
\begin{rmk}
	Here we use the name quantizable function to denote a section of the Weyl bundle, instead of its symbol. The main reason is that, in general, a $D_{BT,k}$-flat section is {\em not} uniquely determined by its symbol, as shown in \cite{ChaLeuLi2023}. The uniqueness only holds when $k>\!\!>1$. However, the name quantizable function will be justified after we show in Corollary \ref{corollary: holomorphic-differential-operator-symbol} that even there are two flat sections with the same symbol (they are different as global sections of sheaf of differential operators on $L^{\otimes k}$), their actions on the Hilbert space $\mathcal{H}_k$ actually coincide. 
\end{rmk}

\begin{rmk}
	The reason why we need the above discussion on formal and non-formal (level $k$) quantizable functions is that in the proof of the orthogonality relations (Theorem \ref{theorem: orthogonality-relations}), we were treating $f$ as a formal smooth function (although there was no $\hbar$) and were considering its associated $D_{BT}$-flat section of $\W_{X,\C}[[\hbar]]$. In particular, only in the formal setting we can define weight components of its associated flat section. So for a non-formal level $k$ quantizable function, we want to think of it as the evaluation of a formal function at 
	$\hbar=\tfrac{\sqrt{-1}}{k}$ (which is always possible locally).
\end{rmk}

Using the compatibility between the flat connections on Weyl bundles and Bargmann--Fock modules, we obtain the following:
\begin{prop}\label{proposition: Bargmann-Fock-differential-operator}
	A level $k$ quantizable function $\alpha$ acts on the Hilbert space $\mathcal{H}_k$ as a holomorphic differential operator. We will denote this differential operator by $P_{\alpha,k}$. 
\end{prop}
\begin{proof}
	Fix $s \in \mathcal{H}_k$. We consider the following Bargmann--Fock action:
	$$
	\alpha\circledast_k O_s.
	$$
	Since the Fedosov connections $D_{BT,k}$ on the Weyl bundle and $D_{BF,k}$ on the Bargmann--Fock sheaf are compatible in the sense that
	\begin{equation}\label{equation: Fedosov-connection-compatibility}
		D_{BF,k}(\alpha\circledast_k O_s)=D_{BT,k}(\alpha)\circledast_k O_s+\alpha\circledast_k D_{BF,k}(O_s)=0,
	\end{equation}
	the output of the action of a quantizable function on a flat section of $\mathcal{F}_{L^{\otimes k}}$ is still flat. Thus by Theorem \ref{theorem: flat-section-Hilbert-space} $\sigma_{L^{\otimes k}}(\alpha \circledast_k O_s)$ is a holomorphic section of $L^{\otimes k}$. This action $s \mapsto \sigma_{L^{\otimes k}}(\alpha \circledast_k O_s)$ is obviously local, since the value of $O_s$ at a point solely depends on the (infinite) jet of $s$ at that point. Thus quantizable functions act as holomorphic differential operators.
\end{proof}

Indeed, we obtain a homomorphism of sheaves of algebras:
\begin{equation}\label{equation: isomorphism-quantizable-function-differential-operator}
	\begin{aligned}
		\varphi: \mathcal{C}_{q,k}^\infty&\rightarrow \mathcal{D}_{L^{\otimes k}}\\
		\alpha &\mapsto P_{\alpha, k}
	\end{aligned}
\end{equation}
One of the main results of \cite{ChaLeuLi2023} is that the map $\varphi$ is actually an isomorphism. 

\begin{thm}[\cite{ChaLeuLi2023}, Theorem 1.3]\label{theorem: almost-holomorphic-function-differential-operators}
	Suppose $X$ is a prequantizable K\"ahler manifold. Then for any positive integer $k$, the homomorphism $\varphi$ in equation \eqref{equation: isomorphism-quantizable-function-differential-operator}
	from the sheaf of algebras of level $k$ quantizable functions to the sheaf of holomorphic differential operators on $L^{\otimes k}$ is an isomorphism. Furthermore, this isomorphism is compatible with the filtration on quantizable functions and that on differential operators by orders, and hence gives an isomorphism of twisted differential operators (TDO).
\end{thm}



The following proposition describes the close relation between the two notions of quantizable functions.
\begin{prop}\label{proposition: quantizable-function-locally-formal}
	If $f\in C^\infty(X)[\hbar]$ is a formal quantizable function on $X$, then the evaluation of $O_f$
	at $\hbar=\tfrac{\sqrt{-1}}{k}$ is a level $k$ (non-formal) quantizable function. On the other hand, if $\alpha$ is any level $k$ (non-formal) quantizable function with symbol $f=\sigma(\alpha)$, then locally on some $U\subset X$, it is the evaluation of the flat section associated to a formal quantizable function $f_0+\hbar\cdot f_1\cdots+ \hbar^l\cdot f_l $ at $\hbar = \tfrac{\sqrt{-1}}{k}$ and there is $f|_U=f_0+\hbar\cdot f_1\cdots+ \hbar^l\cdot f_l $. 
\end{prop}
\begin{proof}
	The first statement is almost obvious from the definition. 
	We now explain that for every level $k$ quantizable function $\alpha$, its symbol $f=\sigma(\alpha)$ is locally the evaluation of a formal quantizable function at 
	$\hbar=\tfrac{\sqrt{-1}}{k}$. By Theorem \ref{theorem: almost-holomorphic-function-differential-operators}, the map $\varphi$ in equation \eqref{equation: isomorphism-quantizable-function-differential-operator} is an isomorphism of sheaves of algebras. We know that the sheaf $\mathcal{D}_{L^{\otimes k}}$ of holomorphic differential operators is locally generated as algebras by holomorphic functions and the operators $\frac{\partial}{\partial z^i}$'s with respect to a local holomorphic frame $e_L^{\otimes k}$ of $L^{\otimes k}$. For a local holomorphic function $g(z)\in\mathcal{O}_X(U)$, it is shown in Appendix \ref{appendix: Fedosov-connection-via-L-infinity-structure} that the corresponding flat section $O_g$ is of weight $0$ and is trivially formal quantizable (see equation \eqref{equation: flat-section-holomorphic-function}). In \cite{ChaLeuLi2023}, we also constructed local function $u^i$'s whose corresponding flat section $O_{u^i}$ acts as $\frac{\partial}{\partial z^i}$ on $L^{\otimes k}$, with respect to a chosen holomorphic frame $e_{L^{\otimes k}}$. It is almost obvious that $u^i$'s are evaluations of formal quantizable functions of order $1$. 
\end{proof}
\begin{rmk}\label{remark: global_obstruction}
	In general, for a global level $k$ quantizable function, it may {\em not} be globally the evaluation of a formal quantizable function at $\hbar=\frac{\sqrt{-1}}{k}$. There are cohomological obstructions to the existence of such a formal quantizable function.  
\end{rmk}

\subsection{Kostant-Souriau operators: old and new}\label{subsection: Kostant-Souriau-operators}
\

To define an associated operator on the Hilbert spaces $\mathcal{H}_k$'s for a general smooth function $f\in C^\infty(X)$, a natural idea is to consider the convergence of the corresponding infinite sum. We give the following example, which says that even for a real analytic function, the corresponding infinite sum might be divergent, as illustrated in the following example:
\begin{eg}\label{example: divergence-infinite-sum}
	We consider the function on $\C$: 
	$$
	f(z,\bar{z})=e^{z\bar{z}}=1+\frac{z\bar{z}}{1!}+\frac{z^2\bar{z}^2}{2!}+\frac{z^3\bar{z}^3}{3!}+\cdots
	$$
	It is easy to see that the infinite sum of the associated Bargmann-Fock action on the function $1$ is
	\begin{align*}
		P_f(1)=&1+\frac{1}{1!}\partial_z(z)+\frac{1}{2!}(\partial_z)^2(z^2)+\frac{1}{3!}(\partial_z)^3(z^3)+\cdots\\
		=&1+1+1+1+\cdots.
	\end{align*}
	which is divergent. 
\end{eg}

Thus it is necessary to construct operators on the Hilbert spaces $\mathcal{H}_k$'s associated to a smooth function $f$ via a truncation of the infinite sum. 
\begin{defn}\label{definition: higher-Kostant-Souriau-operator}
	Given a smooth function $f\in C^\infty(X)$. Then for every $m\geq 0$, by taking the weight $m$ component $\alpha:=(O_f)_m$ of the flat section associated to the function $f$ in Definition \ref{definition: differential-operators-finite-weight-section-Weyl}, we obtain a differential operator: 
	\begin{equation}\label{equation: differential-operator-associated-function}
		P_{f,k,m}:=P_{\alpha,k}: H^0(X, L^{\otimes k})\rightarrow \Gamma(X, L^{\otimes k}). 
	\end{equation}
	The operators $P_{f,k,m}$'s and the partial sums
	$$
	P_{f,k,\leq m}:=\sum_{l=0}^m P_{f,k,l}. 
	$$
	are called {\em higher Kostant--Souriau operators.}
\end{defn}

Here we give some examples of these operators for small $m$:
\begin{eg}
	The operator $P_{f,k,m=0}$ is just the multiplication by $f$. 
\end{eg}

\begin{eg}\label{example: weight-1-differential-operator}
	The weight $1$ component $(O_f)_1$ is explicitly:
	$$
	(O_f)_1=\frac{\partial f}{\partial\bar{z}^j}\bar{y}^j+\frac{\partial^2 f}{\partial z^i\partial\bar{z}^j}y^i\bar{y}^j+\sum_{l\geq 1}(\tilde{\nabla}^{1,0})^l\left(\frac{\partial^2 f}{\partial z^i\partial\bar{z}^j}y^i\bar{y}^j\right).
	$$
	For $s \in \mathcal{H}_k$, we have the following explicit computation:
	\begin{align*}
		P_{f,k,m=1}(s)=&\sigma_{L^{\otimes k}}\left((O_f)_1\circledast_k O_s\right)\\
		=&\sigma_{L^{\otimes k}}\left(\left(\frac{\partial f}{\partial\bar{z}^j}\bar{y}^j+\frac{\partial^2 f}{\partial z^i\partial\bar{z}^j}y^i\bar{y}^j\right)\circledast_k O_s\right)\\
		=&\sigma_{L^{\otimes k}}\left(\left(\frac{\partial f}{\partial\bar{z}^j}\bar{y}^j+\frac{\partial^2 f}{\partial z^i\partial\bar{z}^j}y^i\bar{y}^j\right)\circledast_k \left(s+\tilde{\nabla}^{1,0}(s)\right)\right)\\
		=&\sigma_{L^{\otimes k}}\left( \frac{\partial f}{\partial\bar{z}^j}\bar{y}^j\circledast_k (y^i\otimes\nabla_{\frac{\partial}{\partial z^i}}(s)) \right)+\frac{\sqrt{-1}}{k} \omega^{\bar{j}i}\frac{\partial^2 f}{\partial z^i\partial\bar{z}^j}\cdot s\\
		=&\frac{\sqrt{-1}}{k}\omega^{\bar{j}i}\frac{\partial f}{\partial\bar{z}^j}\cdot\nabla_{\frac{\partial}{\partial z^i}}(s)+\frac{\sqrt{-1}}{k}\omega^{\bar{j}i}\frac{\partial^2 f}{\partial z^i\partial\bar{z}^j}\cdot s\\
		=&\frac{\sqrt{-1}}{k}\nabla_{X_f^{1,0}}(s)+\frac{1}{k}\Delta f\cdot s\\
		=&\frac{\sqrt{-1}}{k}\nabla_{X_f}(s)+\frac{1}{k}\Delta f\cdot s.
\end{align*}
The second equality we have used the following vanishing by the definition of Bargmann-Fock action and the symbol map $\sigma_{L^{\otimes k}}$: for every $l\geq 1$, 
$$
\sigma\left((\tilde{\nabla}^{1,0})^l\left(\frac{\partial^2 f}{\partial z^i\partial\bar{z}^j}y^i\bar{y}^j\right)\circledast_k O_s\right)=0.
$$	
	In the third equality we have used the expression of $O_s$ in Proposition \ref{proposition: expression-flat-section-Fock-bundle}, and also the fact that higher order terms of $O_s$ will vanish after taking the symbol map $\sigma_{L^{\otimes k}}$. Here $X_f^{1,0}=\omega^{\bar{j}i}\frac{\partial f}{\partial\bar{z}^j}\frac{\partial}{\partial z^i}$ denotes the $(1,0)$ component of the Hamiltonian vector field $X_f$ associated to $f$ and $\Delta f = \sqrt{-1}\omega^{\bar{j}i}\frac{\partial^2 f}{\partial z^i\partial\bar{z}^j}$ denotes the Dolbeault Laplacian of $f$. 
\end{eg}

Recall that the {\em Kostant--Souriau pre-quantum operator} associated to $f$ is given by
\begin{equation*}
	H_{f,k} := f + \frac{\sqrt{-1}}{k} \nabla_{X_f}: \Gamma(X, L^{\otimes k}) \to \Gamma(X, L^{\otimes k}).
\end{equation*}
We can see from the above examples that, as operators from $\mathcal{H}_k \to \Gamma(X, L^{\otimes k})$, we have
\begin{equation*}
	P_{f, k, \leq 1} = H_{f,k} + \frac{1}{k} \Delta f \cdot,
\end{equation*}
i.e., $P_{f,k, \leq1}$ only differs from the classical Kostant--Souriau pre-quantum operators $H_{f,k}$ by a correction term $\tfrac{1}{k} \Delta f$.

\begin{eg}\label{example: weight-2-differential-operator}
	By a tedious differential-geometric computation, we obtain the following formula for any $s \in \mathcal{H}_k$:
	\begin{equation*}
		P_{f, k, m=2}(s) = \frac{1}{k^2} \left( -\nabla_{G_f} s + \sqrt{-1} \nabla_{X_{\Delta f}^{1, 0}} s + \frac{1}{2} (\Delta^2 f) \cdot s \right).
	\end{equation*}
	Here, $G_f \in \Gamma(X, \operatorname{Sym}^2 TX)$ is obtained from $(\tilde{\nabla}^{0, 1})^2 f \in \Gamma(X, \operatorname{Sym}^2 \overline{T^*X})$ via the isomorphism $\omega^{-1}: \overline{T^*X} \to TX$, and $\nabla_{G_f}$ denotes the higher covariant derivative introduced in Lemma/Definition \ref{definition: higher-covariant-derivative}. In this sense, $G_f$ may be viewed as a higher-order analogue of the $(1, 0)$-component $X_f^{1, 0}$ of the Hamiltonian vector field associated to $f$. 
\end{eg}

\subsection{Norm estimates for the Kostant-Souriau operators}
\

We first give an estimate for the operators $P_{f,k,m}$ with respect to the $L^2$-norms on sections of the pre-quantum line bundles. 
	\begin{prop}\label{proposition: norm-estimated-differential-operators-f}
		There exists a constant $C_{f,m} > 0$, depending only on $f$ and $m$, such that for any $s\in H^0(X, L^{\otimes k})$
		\begin{equation}\label{equation: norm-estimate-Bargmann-Fock-action}
			\lVert P_{f,k,m}(s)\lVert_{L^2} \leq C_{f,m}\cdot k^{-m/2}\cdot \lVert s\lVert_{L^2}. 
		\end{equation}	
	\end{prop}
	\begin{proof}
		The part of the weight $m$ component of $(O_f)_m$ which contribute non-trivially to $P_{f, k, m} = P_{(O_f)_m,k}$ is the sum of terms 
		$$
		\alpha\in\hbar^p\cdot\Sym^q\overline{T^*X}\otimes\Sym^r T^*X,
		$$	
		where $p+q=m$ and $r\leq q$. By Theorem \ref{theorem: estimate-higher-covariant-derivative} and Proposition \ref{proposition: Bargmann-Fock-action-higher-covariant-derivative}, there exists a constant $C_\alpha$ independent of the quantum level $k$, such that 
		$$
		\lVert P_{\alpha,k} \rVert_{op} \leq C_\alpha\cdot\frac{1}{k^{p+q}}\cdot k^{\frac{q-r}{2}}=C_\alpha\cdot k^{-m+\frac{q-r}{2}}\leq C_\alpha\cdot k^{-\frac{m}{2}}.
		$$
		Equation \eqref{equation: norm-estimate-Bargmann-Fock-action} follows by taking the sum of these inequalities for $\alpha$'s. 
	\end{proof}
	\begin{rmk}
		It is clear from Proposition \ref{proposition: Bargmann-Fock-action-higher-covariant-derivative} that the operators $P_{f,k,m}$'s associated to a function $f$ can naturally be extended to differential operators acting on smooth sections $\Gamma(X, L^{\otimes k})$. However, the norm estimate in Proposition \ref{proposition: norm-estimated-differential-operators-f} holds only when we regard $P_{f,k,m}$ as operators acting on the Hilbert spaces $\mathcal{H}_k = H^0(X, L^{\otimes k})$. 
	\end{rmk}
	
\begin{lem}\label{lemma: partial-bar-D-f-on-s}
	For any smooth function $f\in C^\infty(X)$, holomorphic section $s\in H^0(X, L^{\otimes k})$, and any $m\geq 0$,  we have
	$$
	\bar{\partial}(P_{f,k,\leq m}s)=\sigma_{L^{\otimes k}}\left(\delta^{0,1}(O_f)_{m+1}\circledast_k O_s \right).
	$$
\end{lem}
\begin{proof}
	Recall that $P_{f,k,\leq m}s:=\sigma_{L^{\otimes k}}\left((O_f)_{\leq m}\circledast_k O_s\right)$. First, the following equality 
	$$
	\bar{\partial}\left(P_{f,k,\leq m}(s)\right)=\sigma_{L^{\otimes k}}\circ D_{BF,k}^{0,1}\left((O_f)_{\leq m}\circledast_k O_s\right)
	$$
	follows by substituting $\gamma = (O_f)_{\leq m}\circledast_k O_s$ in the following Lemma \ref{lemma: symbol-cochain-map}. 
	
	Next, in the proof of Lemma \ref{lemma: iteration-relation-weight-components-flat-section}, we have seen the following equation: 
	\begin{equation}\label{equation: 0-1-component-Fedosov-connection}
		D_{BT}^{0,1}(O_f)_{\leq m}=\delta^{0,1} (O_f)_{m+1}.
	\end{equation}	
	By using the compatibility between the Fedosov connection on Bargmann--Fock and Weyl bundles, we get 
	\begin{align*}
		D_{BF,k}^{0,1}\left((O_f)_{\leq m}\circledast_k O_s\right)=&(D_{BT}^{0,1}(O_f)_{\leq m})\circledast_k O_s+(O_f)_{\leq m}\circledast_k D_{BF,k}^{0,1}(O_s)\\
		=&(D_{BT}^{0,1}(O_f)_{\leq m})\circledast_k O_s\\
		=&(\delta^{0,1}(O_f)_{m+1})\circledast_k O_s.
	\end{align*}
	Here, in the second equality, we are using the flatness $D_{BF,k}(O_s)=0$, 
	and the last equality follows from equation \eqref{equation: 0-1-component-Fedosov-connection}. To summarize, by combining the previous computations, we obtain that 
	\begin{equation*}
		\bar{\partial}(P_{f,k,\leq m}s)=\sigma_{L^{\otimes k}} \circ D_{BF,k}^{0,1}\left((O_f)_{\leq m}\circledast_k O_s\right) = \sigma_{L^{\otimes k}} \left(\delta^{0,1}(O_f)_{m+1}\circledast_k O_s \right).
	\end{equation*}
\end{proof}
\begin{rmk}
	In the proof of this lemma, we can see that the flatness of $O_f$ under the Fedosov connection plays an indispensable role.  
\end{rmk}
\begin{lem}\label{lemma: symbol-cochain-map}
	For every section $\gamma$ of the level $k$ Bargmann--Fock sheaf, we have
	$$
	\sigma_{L^{\otimes k}}\circ D_{BF,k}^{0,1}(\gamma)=\bar{\partial}\circ\sigma_{L^{\otimes k}}(\gamma),
	$$
	i.e., the symbol map $\sigma_{L^{\otimes k}}$ is a cochain map with respect to the differentials $D_{BF,k}^{0,1}$ and $\bar{\partial}$.  
\end{lem}
\begin{proof}
	From the expression of the Fedosov connection $D_{BF,k}$, it is clear that $D_{BF,k}^{0,1}$ consists of three parts: the first part is $\bar{\partial}$, the second part is the Bargmann--Fock action $I_m\circledast_k-$, and the last part is $J_{\alpha,m}\circledast_k-$ for $\alpha=Ric_X$. (See Appendix \ref{subsection: Fedosov-quantization} for the notations).  It is easy to see that for any $m\geq 2$, we have the following vanishing by type reasons
	$$
	\sigma_{L^{\otimes k}}\left(I_m\circledast_k\gamma\right)=0,
	$$
	since the degree of $T^*X$ is strictly greater than $\overline{T^*X}$ in $I_m$. Similarly, there are also the vanishing for all $m\geq 1$:
	$$
	\sigma_{L^{\otimes k}}\left(J_{\alpha,m}\circledast_k\gamma\right)=0,
	$$
	and these vanishing results imply the lemma. 
\end{proof}

By combining Lemma \ref{lemma: partial-bar-D-f-on-s} and the estimate in Proposition \ref{proposition: norm-estimated-differential-operators-f}, we immediately obtain the following proposition which says that for any $s\in H^0(X, L^{\otimes k})$ of $L^2$-norm $1$, $\lVert\bar{\partial}(P_{f,k,\leq l}(s))\lVert_{L^2}$ is {\em asymptotically} small as $k\rightarrow\infty$:
\begin{prop}\label{proposition: estimate-partial-bar-D-f-action-s}
	There exists a constant $\tilde{C}_{f,m}$, depending only on the function $f$ and the weight $m$, such that 
	\begin{equation}\label{equation: estimate-partial-bar-D-f-action-s}
		\lVert \bar{\partial}(P_{f,k,\leq m}s) \rVert_{L^2} \leq \tilde{C}_{f,m}\cdot k^{-\frac{m}{2}}\cdot \lVert s \rVert_{L^2}.
	\end{equation}
\end{prop}

\section{Beyond Tuynman's Lemma: orthogonality relations}


In this section, we 
establish a set of orthogonality relations for the higher Kostant--Souriau operators $P_{f,k,m}$'s, generalizing the following classical Tuynman's Lemma.
\begin{lem}[\cite{Tuynman87}]
	For every smooth function $f$, and any two holomorphic sections $s_1,s_2\in \mathcal{H}_k$, we have
	$$
	\left\langle \frac{\sqrt{-1}}{k}\nabla_{X_f}s_1, s_2\right\rangle=\left\langle -\frac{1}{k}\Delta f\cdot s_1, s_2\right\rangle.
	$$
\end{lem}
This implies that $\Pi_k \circ H_{f,k} = T_{f - \frac{1}{k} \Delta f}$.
\begin{rmk}
The original statement of Tuynman's Lemma (\cite{Tuynman87}; see also \cite[Proposition 4.1]{BorHopSchSch1991}) is for level $k=1$; here we are stating the all-level statement. 
\end{rmk}

From the computation in Example \ref{example: weight-1-differential-operator}, it is immediate to see that Tuynman's Lemma can be rephrased as the following orthognality relation: for any $s_1,s_2\in\mathcal{H}_k$, and for any smooth function $f\in C^\infty(X)$, we have
\begin{equation}\label{equation: Tuyman-lemma-refrased}
\langle P_{f,k,m=1}(s_1), s_2\rangle=0.
\end{equation}
The main result in this section is that Tuynman's Lemma can be generalized to the following series of orthogonality relations:
\begin{equation}\label{equation: identification-differential-operators}
\langle P_{f,k,m}(s_1), s_2\rangle=0,\quad m\geq 1.
\end{equation}
\begin{rmk}
 Interestingly, although the above statement holds for every fixed level $k$ (which is geometrically the tensor power of the pre-quantum line bundle), we need to consider the formal deformation quantization because only from there we can see the weight structures. 
\end{rmk}

\subsection{Another series of differential operators}
\

In the previous subsection, we defined the higher Kostant--Souriau differential operators $P_{f,k,m}$, $m\geq 0$ associated to a function $f\in C^\infty(X)$ (Definition \ref{definition: higher-Kostant-Souriau-operator}). Here, for each $m \geq 1$, we define another differential operator 
$$\tilde{P}_{f,k,m}:\mathcal{H}_k\rightarrow \Gamma(X, L^{\otimes k})$$ 
via the following equation: for any $s\in\mathcal{H}_k$, 
\begin{equation}\label{equation: differential-operator-P-tilde}
	\tilde{P}_{f,k,m}(s)\cdot\omega^n=\nabla^{1,0}\circ\sigma_{L^{\otimes k}} \left(\delta^{0,1}(\hbar O_f)_{m+1}\circledast_k O_s\right)\cdot n\omega^{n-1}.
\end{equation}

In fact, these two series of differential operators are identical:
\begin{prop}\label{proposition: gamma-m-differential}
	Let $f\in C^\infty(X)$ be a smooth function on $X$, and $s\in H^0(X, L^{\otimes k})$ be a holomorphic section of $L^{\otimes k}$. Then for every weight $m\geq 1$, we have the following equality:
	$$
	P_{f,k,m}(s)=\tilde{P}_{f,k,m}(s). 
	$$
\end{prop}

The proof of Proposition \ref{proposition: gamma-m-differential} consists of two parts. In the first part, we show that on the flat space $\mathbb{C}^n$ with the standard K\"ahler structure, the statement of Proposition \ref{proposition: gamma-m-differential} holds. In the second part, we show that on a compact K\"ahler manifold, the proof can be reduced to the flat space situation. 

\subsubsection{A technical lemma}
\

	Let $\alpha$ be a section of the Weyl bundle on $\mathbb{C}^n$ such that it is a polynomial in $\bar{y}^j$'s and $\hbar$. Let $\omega=\sqrt{-1}\sum_i dz^i\wedge d\bar{z}^i$ denote the standard K\"ahler form on $\mathbb{C}^n$. 
	
	For such an $\alpha$, we will define two operators acting on the Bargmann--Fock sheaf on $\mathbb{C}^n$. The first one is simply the Bargmann--Fock action 
	$$
	T_{\alpha,k}(\beta):=\alpha\circledast_k \beta.
	$$
	The second operator $\tilde{T}_{\alpha,k}$ is defined via the following equality:
	$$
	\tilde{T}_{\alpha,k}(\beta)\cdot\omega^n= d\left(\delta^{0,1}(\hbar \alpha)\circledast_k \beta\right)\cdot (n\omega^{n-1}) = \partial\left(\delta^{0,1}(\hbar \alpha)\circledast_k \beta\right)\cdot (n\omega^{n-1}).
	$$
	Here $d$ denotes the de Rham differential on $X$, which is also extended to the trivialized Weyl bundle. 
\begin{lem}\label{lemma: technical-lemma-flat-space}
Suppose that $\alpha$ and $\beta$ satisfy the following equations respectively:
\begin{equation}\label{equation: alpha-holomorphic-Taylor-expansion}
(\partial-\delta^{1,0})\alpha=0, \hspace{6mm}(\partial-\delta^{1,0})\beta=0. 
\end{equation}
Suppose furthermore that every term of $\alpha$ has at least polynomial degree $1$ in $\overline{y}^j$'s. There is the following equality:
$$
T_{\alpha,k}(\beta)=\tilde{T}_{\alpha,k}(\beta).
$$
\end{lem}
\begin{proof}
Without loss of generality, we assume that $\alpha$ does not involve the formal variable $\hbar$ and is a homogeneous polynomial in $\overline{y}^j$'s of degree $l > 0$. We write $\alpha, \beta$ explicitly as
\begin{align*}
	\alpha=&\sum_{p\geq 0} \alpha_{i_1 \cdots i_p, \bar{j}_1 \cdots \bar{j}_l} y^{i_1}\cdots y^{i_p}\bar{y}^{j_1}\cdots\bar{y}^{j_l},\\
	\beta=&\sum_{q\geq 0} \beta_{k_1\cdots k_q}y^{k_1}\cdots y^{k_q}.
\end{align*}
We first explain how equation \eqref{equation: alpha-holomorphic-Taylor-expansion} is reflected on the indices $\alpha_{i_1\cdots i_p, \bar{j}_1\cdots \bar{j}_l}$ and $\beta_{k_1\cdots k_q}$'s. It is easy to see that equation \eqref{equation: alpha-holomorphic-Taylor-expansion} implies that
$$
\frac{\partial\alpha_{i_1\cdots i_p, \bar{j}_1\cdots \bar{j}_l}}{\partial z^r}=(p+1)\cdot\alpha_{i_1\cdots i_p r, \bar{j}_1\cdots \bar{j}_l}.
$$
Here we have a detailed computation of the coefficient $\tfrac{1}{p+1}$: from equation \eqref{equation: alpha-holomorphic-Taylor-expansion}, we have 
$$
\partial(\alpha_{i_1\cdots i_p, \bar{j}_1\cdots \bar{j}_l}y^{i_1}\cdots y^{i_p}\bar{y}^{j_1}\cdots\bar{y}^{j_l})=\delta^{1,0}(\alpha_{i_1\cdots i_pi_{p+1}, \bar{j}_1\cdots \bar{j}_l}y^{i_1}\cdots y^{i_p}y^{i_{p+1}}\bar{y}^{j_1}\cdots\bar{y}^{j_l}),
$$
which implies that 
$$
\frac{\partial\alpha_{i_1\cdots i_p, \bar{j}_1\cdots \bar{j}_l}}{\partial z^{i_{p+1}}}\cdot dz^{i_{p+1}}\otimes y^{i_1}\cdots y^{i_p}\bar{y}^{j_1}\cdots\bar{y}^{j_l}=(p+1)\cdot\alpha_{i_1\cdots i_pi_{p+1}, \bar{j}_1\cdots \bar{j}_l}dz^{i_{p+1}}\otimes y^{i_1}\cdots y^{i_p}\bar{y}^{j_1}\cdots\bar{y}^{j_l}.
$$
Similarly, we have 
$$
\frac{\partial \beta_{k_1\cdots k_q}}{\partial z^r}=(q+1)\cdot \beta_{k_1\cdots k_qr}.
$$

The Bargmann--Fock action of $\alpha$ is via the following fiberwise differential operator
$$
\sum_{p \geq 0} \frac{1}{k^l} \alpha_{i_1\cdots i_p, \bar{j}_1\cdots \bar{j}_l} \cdot \frac{\partial^l}{\partial y^{j_1} \cdots \partial y^{j_l}}\circ m_{y^{i_1} \cdots y^{i_p}}. 
$$
And that of $\delta^{0,1}\alpha$ is explicitly given as
$$
\sum_{p\geq 0} \frac{1}{k^{l-1}} \cdot l \cdot \alpha_{i_1\cdots i_p, \bar{j}_1\cdots \bar{j}_l} d\bar{z}^{j_l}\cdot \frac{\partial^{l-1}}{\partial y^{j_1} \cdots \partial y^{j_{l-1}}}\circ m_{y^{i_1} \cdots y^{i_p}}. 
$$
(Here we are implicitly using the symmetry of the indices $\bar{J}$). 
There is the following tedious but straightforward computation:
\begin{align*}
	&\tilde{T}_{\alpha,k}(\beta) \cdot \omega^n\\
	=&\partial\left(\sum_{p,q\geq 0}\frac{1}{k^l} \cdot l \cdot \alpha_{i_1\cdots i_p, \bar{j}_1\cdots \bar{j}_l} \beta_{k_1\cdots k_q} d\bar{z}^{j_l}\cdot \frac{\partial^{l-1} \left( y^{i_1} \cdots y^{i_p} y^{k_1}\cdots y^{k_q} \right)}{\partial y^{j_1} \cdots \partial y^{j_{l-1}}} \right) \cdot n\sqrt{-1}\omega^{n-1}\\
		=& \sum_{p,q\geq 0}\frac{1}{k^l} \cdot l \cdot \frac{\partial\alpha_{i_1\cdots i_p, \bar{j}_1\cdots \bar{j}_l}}{\partial z^r} \beta_{k_1\cdots k_q} \frac{\partial^{l-1} \left( y^{i_1} \cdots y^{i_p} y^{k_1}\cdots y^{k_q} \right)}{\partial y^{j_1} \cdots \partial y^{j_{l-1}}} \cdot dz^r \wedge d\bar{z}^{j_l} \wedge n\sqrt{-1}\omega^{n-1}\\
		&+ \sum_{p,q\geq 0}\frac{1}{k^l} \cdot l \cdot \alpha_{i_1\cdots i_p, \bar{j}_1\cdots \bar{j}_l} \frac{\partial \beta_{k_1\cdots k_q}}{\partial z^r} \frac{\partial^{l-1} \left( y^{i_1} \cdots y^{i_p} y^{k_1}\cdots y^{k_q} \right)}{\partial y^{j_1} \cdots \partial y^{j_{l-1}}} \cdot dz^r \wedge d\bar{z}^{j_l} \wedge n\sqrt{-1}\omega^{n-1}\\
		=&\sum_{p,q\geq 0}\frac{1}{k^l} \cdot l \cdot \frac{\partial\alpha_{i_1\cdots i_p, \bar{j}_1\cdots \bar{j}_l}}{\partial z^{j_l}} \beta_{k_1\cdots k_q} \frac{\partial^{l-1} \left( y^{i_1} \cdots y^{i_p} y^{k_1}\cdots y^{k_q} \right)}{\partial y^{j_1} \cdots \partial y^{j_{l-1}}} \cdot \omega^n\\
		&+\sum_{p,q\geq 0}\frac{1}{k^l} \cdot l \cdot \alpha_{i_1\cdots i_p, \bar{j}_1\cdots \bar{j}_l} \frac{\partial \beta_{k_1\cdots k_q}}{\partial z^{j_l}} \frac{\partial^{l-1} \left( y^{i_1} \cdots y^{i_p} y^{k_1}\cdots y^{k_q} \right)}{\partial y^{j_1} \cdots \partial y^{j_{l-1}}} \cdot \omega^n.
\end{align*}
In the above computation, we have used the following equations:
\begin{align*}
	\omega^n&=n!\cdot (\sqrt{-1})^n \cdot dz^1\wedge d\bar{z}^1\wedge\cdots\wedge dz^n\wedge d\bar{z}^n;\\
	dz^i\wedge d\bar{z}^j \wedge \omega^{n-1} &=0, \text{ for } i \neq j;\\
	dz^j\wedge d\bar{z}^j \wedge \omega^{n-1} & = (n-1)! \cdot (\sqrt{-1})^{n-1} \cdot dz^1\wedge d\bar{z}^1\wedge\cdots\wedge dz^n\wedge d\bar{z}^n=\frac{1}{n \sqrt{-1}}\cdot\omega^n.
\end{align*}
On the other hand, we have 
\begin{align*}
	&T_{\alpha,k}(\beta)\\
	=& \sum_{p,q\geq 0}\frac{1}{k^l} \cdot \alpha_{i_1\cdots i_p, \bar{j}_1\cdots \bar{j}_l} \beta_{k_1\cdots k_q} \frac{\partial^l}{\partial y^{j_1} \cdots \partial y^{j_l}} \left( y^{i_1} \cdots y^{i_p} y^{k_1} \cdots y^{k_q} \right)\\
	=&\sum_{p>0,q\geq 0}\frac{1}{k^l} \cdot \alpha_{i_1\cdots i_p, \bar{j}_1\cdots \bar{j}_l}\beta_{k_1\cdots k_q} \cdot l\cdot p\cdot\frac{\partial y^{i_p}}{\partial y^{j_l}} \cdot\frac{\partial^{l-1}}{\partial y^{j_1} \cdots \partial y^{j_{l-1}}} \left( y^{i_1} \cdots y^{i_{p-1}} y^{k_1} \cdots y^{k_q} \right)\\
	&+\sum_{p\geq 0,q>0}\frac{1}{k^l} \cdot \alpha_{i_1\cdots i_p, \bar{j}_1\cdots \bar{j}_l} \beta_{k_1\cdots k_q} \cdot l\cdot q\cdot\frac{\partial y^{k_q}}{\partial y^{j_l}} \cdot \frac{\partial^{l-1}}{\partial y^{j_1} \cdots \partial y^{j_{l-1}}} \left( y^{i_1} \cdots y^{i_p} y^{k_1}\cdots y^{k_{q-1}}\right)\\
	=&\sum_{p>0,q\geq 0}\frac{1}{k^l} \cdot l\cdot(p \alpha_{i_1\cdots i_{p-1}j_l,\bar{j}_1\cdots \bar{j}_l} ) \cdot \beta_{k_1\cdots k_q} \cdot\frac{\partial^{l-1}}{\partial y^{j_1} \cdots \partial y^{j_{l-1}}} \left( y^{i_1} \cdots y^{i_{p-1}} y^{k_1} \cdots y^{k_q} \right)\\
	&+\sum_{p\geq 0,q>0} \frac{1}{k^l} \cdot l \cdot \alpha_{i_1\cdots i_p, \bar{j}_1\cdots \bar{j}_l} \cdot (q \beta_{k_1\cdots k_{q-1}j_l}) \cdot\frac{\partial^{l-1}}{\partial y^{j_1} \cdots \partial y^{j_{l-1}}} \left(y^{i_1} \cdots y^{i_p} y^{k_1}\cdots y^{k_{q-1}}\right)\\
		=&\sum_{p>0,q\geq 0}\frac{1}{k^l} \cdot l \cdot \frac{\partial\alpha_{i_1\cdots i_{p-1},\bar{j}_1\cdots \bar{j}_l}}{\partial  z^{j_l}} \beta_{k_1\cdots k_q} \cdot \frac{\partial^{l-1}}{\partial y^{j_1} \cdots \partial y^{j_{l-1}}} \left( y^{i_1} \cdots y^{i_{p-1}} y^{k_1} \cdots y^{k_q} \right)\\
	&+\sum_{p\geq 0,q>0} \frac{1}{k^l} \cdot l \cdot \alpha_{i_1\cdots i_p, \bar{j}_1\cdots \bar{j}_l} \frac{\partial \beta_{k_1\cdots k_{q-1}}}{\partial z^{j_l}} \cdot\frac{\partial^{l-1}}{\partial y^{j_1} \cdots \partial y^{j_{l-1}}} \left( y^{i_1} \cdots y^{i_p} y^{k_1}\cdots y^{k_{q-1}}\right)\\
	=& \tilde{T}_{\alpha,k}(\beta). 
\end{align*}
In the second step, we have used the condition that $l > 0$ to exclude the case when $p=0$ in the first summation and the case when $q=0$ in the second summation.
\end{proof}


\subsubsection{Proof of Proposition \ref{proposition: gamma-m-differential}}
\

We prove this proposition by showing that for any $s\in\mathcal{H}_k$ and any point $x_0\in X$, we have
$$
P_{f,k,m}(s) \vert_{x_0}=\tilde{P}_{f,k,m}(s) \vert_{x_0}.
$$
We first assume that the K\"ahler form is real analytic, and later we will explain how this technical condition can be removed. Under this technical assumption, we can choose a special holomorphic frame $e_L$ and local coordinates centered at $x_0$, such that we have the vanishing of the coefficients of the Taylor expansion of $\rho=\log h(e_L, e_L)$ at $x_0$:
\begin{equation}\label{equation: local-technical-asumption}
	\begin{aligned}
		\frac{\partial^{\lvert I \rvert}\rho}{\partial z^I}(x_0)=&\frac{\partial^{\lvert J \rvert}\rho}{\partial \bar{z}^J}(x_0)=0,\hspace{3mm}\text{all}\ I,J;\\
		\frac{\partial^{\lvert I \rvert+1}\rho}{\partial\bar{z}^j\partial z^I}(x_0)=&\frac{\partial^{\lvert J \rvert+1}\rho}{\partial z^i\partial \bar{z}^J}(x_0)=0, \hspace{3mm} \lvert I \rvert,\lvert J \rvert\neq 1. 
	\end{aligned}
\end{equation}

We can further require that the K\"ahler form at $x_0$ is the standard one by doing a linear transformation of the coordinate functions:
$$
\omega \vert_{x_0}=\sqrt{-1}\sum_{i=1}^ndz^i\wedge d\bar{z}^i. 
$$
This is called a K-coordinate of $X$ centered at $x_0$. Recall that the pre-quantum condition implies that $\rho$ is a K\"ahler potential, i.e., 
$$
\omega_{i\bar{j}}=\frac{\partial^2\rho}{\partial z^i\partial\bar{z}^j}.
$$
Also recall that the Christoffel symbols of the Levi-Civita connection are given by the following explicit formulas:
$$
\nabla(\partial_{z^j})=\Gamma_{ij}^k\cdot dz^i\otimes\partial_{z^k}=\omega^{\bar{l}k}\cdot\frac{\partial\omega_{j\bar{l}}}{\partial z^i}\cdot dz^i\otimes\partial_{z^k}.
$$
It is then obvious that equation \eqref{equation: local-technical-asumption} implies the vanishing of the Christoffel symbols $\Gamma_{ij}^k$ and all their holomorphic derivatives at the base point $x_0$. 

We can now apply this vanishing to look at the flat section associated to a smooth function $f$ at the base point $x_0$. Let us write $\tilde{\alpha}:=(O_f)_m$ as
$$
\tilde{\alpha}=\sum_{0\leq l\leq m} \hbar^{m-l}\cdot\tilde{\alpha}_{i_1\cdots i_p,j_1,\cdots j_l}y^{i_1}\cdots y^{i_p}\bar{y}^{j_1}\cdots\bar{y}^{j_l}.
$$
The equality $(\nabla^{1,0}-\delta^{1,0})\tilde{\alpha}=0$ then implies the following:
\begin{equation}\label{equation: iteration-equation-O-f}
	\begin{aligned}
		&p\cdot\tilde{\alpha}_{i_1\cdots i_p,j_1,\cdots j_l}dz^{i_p}\otimes y^{i_1}\cdots  y^{i_{p-1}}\cdot\bar{y}^{j_1}\cdots\bar{y}^{j_l} \vert_{x_0}\\
		=&\frac{\partial \tilde{\alpha}_{i_1\cdots i_{p-1},j_1,\cdots j_l}}{\partial z^{i_p}}dz^{i_p}\otimes y^{i_1}\cdots y^{i_{p-1}}\bar{y}^{j_1}\cdots\bar{y}^{j_l} \vert_{x_0}.
	\end{aligned}
\end{equation}
Since $f\in C^\infty(X)$, there is $\sigma(O_f)=\sigma((O_f)_0)=f$. By the hypothesis that $m \geq 1$, we see that $\sigma(\tilde{\alpha}) = \sigma((O_f)_m)=0$. From the equality $(\nabla^{1,0}-\delta^{1,0})\tilde{\alpha}=0$ again, we can deduce that $\tilde{\alpha}_{i_1\cdots i_p, \emptyset} = 0$ for all $p \geq 0$ and all indices $i_1, ..., i_p$. Thus,
$$
\tilde{\alpha}=\sum_{1\leq l\leq m} \hbar^{m-l}\cdot\tilde{\alpha}_{i_1\cdots i_p,j_1,\cdots j_l}y^{i_1}\cdots y^{i_p}\bar{y}^{j_1}\cdots\bar{y}^{j_l}.
$$

Similarly, given $s\in\mathcal{H}_k$, locally we can identify it with a holomorphic function with respect to the frame $e_L^{\otimes k}$, i.e., $s=g(z)\cdot e_L^{\otimes k}$. By Proposition \ref{proposition: expression-flat-section-Fock-bundle}, its associated flat section is locally of the explicit form
$$
\tilde{\beta}:=O_s=\sum_{i\geq 0}(\tilde{\nabla}^{1,0})^i(s)=\sum_{q\geq 0}\tilde{\beta}_{k_1\cdots k_q}\cdot y^{k_1}\cdots y^{k_q}\otimes e_L^{\otimes k}.
$$
A similar argument shows that the coefficients $\tilde{\beta}_{k_1\cdots k_q}$ satisfy the following equation (here we also need the condition on the frame $e_L$):
\begin{equation}\label{equation: iteration-equation-O-s}
q\cdot \tilde{\beta}_{k_1\cdots k_q}dz^{k_q}\otimes y^{k_1}\cdots y^{k_{q-1}} \vert_{x_0}=\frac{\partial\tilde{\beta}_{k_1\cdots k_{q-1}}}{\partial z^{k_q}}dz^{k_q}\otimes y^{k_1}\cdots y^{k_{q-1}} \vert_{x_0}. 
\end{equation}

On the other hand, we can identify a neighborhood $x_0\in U_{x_0}$ with a disk $U\subset \C^n$ centered at the origin using a K-coordinate centered at $x_0$. To apply Lemma \ref{lemma: technical-lemma-flat-space}, we will also define $\alpha$ and $\beta$ as
\begin{align*}
	\alpha := &\sum_{1\leq l \leq m}\hbar^{m-l}\cdot\sum_{i\geq 0}((\delta^{1,0})^{-1}\circ\partial)^i\left(\tilde{\alpha}_{\emptyset,j_1\cdots j_l}\bar{y}^{j_1}\cdots\bar{y}^{j_l}\right),\\
	\beta := &\sum_{i\geq 0}((\delta^{1,0})^{-1}\circ\partial)^i(g(z)).
\end{align*}
From this construction, it is easy to see that $\alpha$ and $\beta$ satisfy equation \eqref{equation: alpha-holomorphic-Taylor-expansion}. It follows from equations \eqref{equation: iteration-equation-O-f} and \eqref{equation: iteration-equation-O-s} that their restrictions to $x_0$ are identical with respect to the K-coordinate and frame $e_L^{\otimes k}$, namely, 
$$
\alpha \vert_{x_0}=\tilde{\alpha} \vert_{x_0};\hspace{5mm}\beta \vert_{x_0}=
(\tilde{\beta}\otimes e_L^{\otimes k}) \vert_{x_0}. 
$$
(Roughly speaking, the covariant derivative and ordinary derivative are identified at the base point $x_0$, by the vanishing of the Christoffel symbols and their holomorphic derivatives at $x_0$.)
\begin{rmk}
There are then the Weyl bundle and Bargmann--Fock sheaf on $U$, either using the pull-back metric from $U_{x_0}$ or the standard K\"ahler metric on $\mathbb{C}^n$. 
\end{rmk}

Also the fiberwise Bargmann--Fock action at $x_0$ are identical with respect to the two K\"ahler metrics. It follows from this observation and Lemma \ref{lemma: technical-lemma-flat-space} that
$$
\tilde{P}_{f,k,m}(s) \vert_{x_0}=\left(\sigma\circ\tilde{T}_{\alpha,k}(\beta)\right)\otimes e_L^{\otimes k} \vert_{x_0} =\left(\sigma\circ T_{\alpha,k}(\beta)\right)\otimes e_L^{\otimes k} \vert_{x_0}
=P_{f,k,m}(s) \vert_{x_0}.$$




To finish the proof, we explain how the technical assumption of the existence of local frame $e_L$ and the K-coordinates centered at $x_0$ can be weakened: since we are fixing the weight $m$ here, we only need to find a frame $e_L$ and a coordinate system at $x_0$, such that the vanishing in equation \eqref{equation: local-technical-asumption} is valid for finitely many $I, J$'s, since the higher order Taylor expansions will be invisible after taking the symbol.
\qed

\subsection{The orthogonality relations}
\

We are now ready to state and prove the main theorem of this section.
\begin{thm}\label{theorem: orthogonality-relations}
Let $f\in C^\infty(X)$. For every weight $m\geq 1$ and 
any holomorphic section $s\in\mathcal{H}_k$,
\begin{equation}\label{equation: Bargmann-Fock-action-lives-in-orthogonal-complement}
	P_{f,k,m}(s)\in\mathcal{H}_k^\perp.
\end{equation}
\end{thm}

\begin{proof}
Let $h$ be the hermitian metric on $L^{\otimes k}$. Consider any $s' \in \mathcal{H}_k$. By Proposition \ref{proposition: gamma-m-differential},
\begin{align*}
	\int_X h(P_{f,k,m}(s), s') \cdot\omega^n=&\int_X h(\tilde{P}_{f,k,m}(s), s') \cdot\omega^n\\
	=&\int_X h(\nabla^{1,0}\circ\sigma_{L^{\otimes k}} \left(\delta^{0,1}(\hbar O_f)_{ m+1}\circledast_k O_s\right),s') \cdot n\omega^{n-1}\\
	=&\int_X\partial_X h(\sigma_{L^{\otimes k}} \left(\delta^{0,1}(\hbar O_f)_{ m+1}\circledast_k O_s\right),s') \cdot n\omega^{n-1}\\
	=&\int_Xd_X\left(h(\sigma_{L^{\otimes k}} \left(\delta^{0,1}(\hbar O_f)_{ m+1}\circledast_k O_s\right),s') \cdot n\omega^{n-1}\right)\\
	=&0.
\end{align*}
Here the third equality follows from the fact that $\bar{\partial}(s')=0$ and the property of Chern connection, and $d_X=\partial_X+\bar{\partial}_X$ denotes the de Rham differential on $X$ and its type decomposition. The last step is saying that 
the integrand $h(P_{f,k,m}(s), s') \cdot\omega^n$
is exact, and the vanishing follows from Stoke's Theorem. 
\end{proof}




\section{Asymptotic holomorphicity of higher Kostant-Souriau operators}\label{section: asymptotic-holomorphicity}


We have seen that the higher Kostant-Souriau operators $P_{f,k,\leq m}$ maps a general holomorphic section $s\in H^0(X, L^{\otimes k})$ to a smooth section of $L^{\otimes k}$. In Proposition \ref{proposition: estimate-partial-bar-D-f-action-s}, we proved that the $L^2$ norm of $\bar{\partial}\circ P_{f,k,\leq m}(s)$ is asymptotically small (comparing to the $L^2$ norm of $s$) which shows that $P_{f,k,\leq m}(s)$ is close to being holomorphic. 

However, this is not adequate from the quantization point of view. A desired {\em quantum operator} should preserve each Hilbert space $\mathcal{H}_k$. We will show in this section that for any $s\in\mathcal{H}_k$, there is an upper bound for the distance between $P_{f,k,\leq m}(s)$ and the subspace $\mathcal{H}_k\subset\Gamma(X, L^{\otimes k})$:
\begin{equation}\label{equation: distance-subspace-holomorphic-section}
\text{dist}(P_{f,k,\leq m}(s), \mathcal{H}_k)=\inf\big\{||P_{f,k,\leq m}(s)-v||_{L^2}: v\in\mathcal{H}_k\big\}\leq C_{f,m}\cdot k^{-\frac{m+1}{2}}\cdot ||s||_{L^2},
\end{equation}
i.e., $P_{f,k,\leq m}$ maps $\mathcal{H}_k$ asymptotically to itself as $k\rightarrow\infty$. Here $C_{f,m}$ is a constant independent of the level $k$ and the holomorphic section $s\in\mathcal{H}_k$. 


This main result is obtained via a comparision between the Kostant-Souriau operators $P_{f,k,\leq m}$ and the Berezin--Toeplitz operators $T_{f,k}$ associated to a same function $f\in C^\infty(X)$ in the following theorem. (Recall that $T_{f,k}(s)$ is defined to be a holomorphic section of $L^{\otimes k}$ and we can choose $v=T_{f,k}(s)$ in equation \eqref{equation: distance-subspace-holomorphic-section} for the estimate).  
\begin{thm}\label{theorem: main}
	For any smooth function $f\in C^\infty(X)$ and any weight $m\geq0$, there exists a constant $C_{f,m} > 0$, depending only on $f$ and $m$, such that for all positive integer $k$, we have
$$
\lVert P_{f,k,\leq m} - T_{f,k} \rVert_{op}\leq C_{f,m} \cdot \frac{1}{k^{\frac{m+1}{2}}}.
$$
In other words, the sequence of differential operators $P_{f,k,\leq m}$ is asymptotic to the Berezin--Toeplitz operator $T_{f,k}$ as $k\rightarrow\infty$. 
\end{thm}

\begin{rmk}
	In Theorem \ref{theorem: main}, $\lVert \cdot \rVert_{op}$ denotes the operator norm of operators $\mathcal{H}_k \to \Gamma(X, L^{\otimes k})$. Although the domain of the operators $P_{f,k,\leq m}-T_{f,k}: \mathcal{H}_k \rightarrow\Gamma(X, L^{\otimes k})$ can be extended to $\Gamma(X, L^{\otimes k})$, their norm estimates only hold when acting on holomorphic sections of $L^{\otimes k}$. 
\end{rmk}

Furthermore, Theorem \ref{theorem: main} provides a more fundamental explanation of why the right hand side of the asymptotic formula  for the composition of Berezin-Toeplitz operators (equation \eqref{equation: BT-star-product-asymptotic}) is given by a sequence of bi-differential operators.

We have seen in Section \ref{subsection: quantizable-function} that if $\alpha\in\Gamma(X, \W_X\otimes\Sym\overline{T^*X})$ is a {\em level $k$ quantizable function}, then it induces a holomorphic differential operator $P_{\alpha,k}$ on the level $k$ Hilbert space $\mathcal{H}_k = H^0(X, L^{\otimes k})$. From Theorem \ref{theorem: main}, it is natural to expect that in this situation, the Berezin-Toeplitz operators and Kostant-Souriau operators (via Bargmann-Fock action) are identical. 
We will prove this in Theorem \ref{theorem: Toeplitz-holomorphic-differential-operator} and describe a characterization of when Berezin--Toeplitz operators $T_{f,k}$ act on the Hilbert space $\mathcal{H}_k$ as holomorphic differential operators.


\subsection{The $m = 1$ case as a consequence of the work of Charles--Polterovich \cite{ChaPol2017}}
\

In \cite{ChaPol2017}, Charles and Polterovich proved a result about the norm estimate of the commutator
\begin{equation*}
	[H_{f,k}, \Pi_k],
\end{equation*}
where $H_{f,k}$ denotes the Kostant--Souriau pre-quantum operator acting on smooth sections of $L^{\otimes k}$. 
In this subsection, we will show that Charles--Polterovich's result implies the $m=1$ case of Theorem \ref{theorem: main}. To see this, we need to reformulate our main result.

\begin{prop}\label{proposition: commutator-Bergman-projection-differential-operators}
	Theorem \ref{theorem: main} is equivalent to the following: 
	there exists $k_0 \in \mathbb{Z}^+$ such that for any $m \geq 0$ and any smooth function $f \in C^\infty(X)$, there exists a constant $C_{f,m} > 0$ (depending only on $f$ and $m$) so that for $k \geq k_0$, we have
	\color{black}
	\begin{equation}\label{equation: main-result-commutator}
		\left\lVert[P_{f,k,\leq m}, \Pi_k]\right\rVert_{op}\leq C_{f,m}\cdot \frac{1}{k^{(m+1)/2}}.
	\end{equation}
\end{prop}
\begin{proof}
	For any holomorphic section $s\in\mathcal{H}_k$, we have
	\begin{align*}
		[P_{f,k,\leq m}, \Pi_k](s)=& P_{f,k,\leq m}\circ\Pi_k(s) - \Pi_k\circ \left(f+\sum_{l=1}^m P_{f,k,l}\right)(s)\\
		=& P_{f,k,\leq m}(s) - \Pi_k (fs)\\
		=&P_{f,k,\leq m}(s) - T_{f,k}(s).
	\end{align*}
	Here we have used the orthogonality relations (Theorem \ref{theorem: orthogonality-relations}) which imply the vanishing $\Pi_k\circ P_{f,k,l}=0$ for all $l>0$ and for all quantum level $k$. 
\end{proof}
Using this proposition, it is easy to see that for $m=1$, the inequality \eqref{equation: main-result-commutator} is a consequence of the following proposition.
\begin{prop}[\cite{ChaPol2017}, Theorem 3.5]\label{proposition: Charles-Polterovich-result}
	There exists a constant $C>0$ such that for any $f\in C^\infty(X)$ and any $k\in\mathbb{Z}^+$,
	\begin{equation}\label{equation: Charles-Polterovich-result}
		\left\lVert [H_{f,k}, \Pi_k] \right\rVert_{op}\leq C\cdot k^{-1} \lvert f \rvert_2,
	\end{equation}
	where $\rvert f \rvert_2$ is the $C^2$-norm of $f$.
\end{prop}

\begin{rmk}
	Indeed, according to the original statement of \cite[Theorem 3.5]{ChaPol2017}, the inequality \eqref{equation: Charles-Polterovich-result} still holds even when $f$ is of class $C^2$ and $\lVert \cdot \rVert_{op}$ is replaced by the operator norm of operators $\Gamma(X, L^{\otimes k}) \to \Gamma(X, L^{\otimes k})$.
\end{rmk}

Let us explain how to prove the $m=1$ case of Theorem \ref{theorem: main} by Proposition \ref{proposition: Charles-Polterovich-result}. In Subsection \ref{subsection: Kostant-Souriau-operators}, we observed that for any $s \in \mathcal{H}_k$, we have
$$
P_{f,k,\leq 1}s = f\cdot s+\frac{\sqrt{-1}}{k}\nabla_{X_f}+\frac{1}{k}\Delta f\cdot s = H_{f,k}(s)+\frac{1}{k} \Delta f \cdot s.
$$
We also observe that
$$
\lVert [\Delta(f), \Pi_k] (s) \rVert = \lVert \Delta(f) \cdot s - T_{\Delta(f), k} s \rVert \leq \lVert \Delta(f) \cdot s \rVert + \lVert T_{\Delta(f), k} s \rVert \leq C_1 \cdot \vert f \rvert_2 \cdot \lVert s \rVert,
$$
for some constant $C_1 > 0$. To conclude, the inequality \eqref{equation: Charles-Polterovich-result}, together with the above observations, imply the inequality in Proposition \ref{proposition: commutator-Bergman-projection-differential-operators} for $m = 1$.

\subsection{Proof of Theorem \ref{theorem: main}}
\

To prove Theorem \ref{theorem: main}, we need a version of the Kodaira--H\"ormander $L^2$ estimate:
\begin{thm}[\cite{Cha2015}, Theorem 2.2]\label{theorem: Hormander-estimate}
	There exist $k_0 \in \mathbb{Z}^+$ and $C > 0$ such that for all $k \in \mathbb{Z}^+$ with $k \geq k_0$, for all $s \in \Gamma(X, L^{\otimes k})$ which is orthogonal to $H^0(X, L^{\otimes k})$, we have
	\begin{equation}
		\lVert s \rVert \leq C\cdot k^{-\frac{1}{2}}\cdot \lVert \overline{\partial} s \rVert,
	\end{equation}
	where $\lVert \overline{\partial} s \rVert$ denotes the natural $L^2$ norm on $\A_X^{0,1}(L^{\otimes k})$ induced by the hermitian metric $h$ on $L^{\otimes k}$ and the K\"ahler metric on $X$. 
\end{thm}

\begin{proof}[Proof of Theorem \ref{theorem: main}]
For any holomorphic section $s \in H^0(X, L^{\otimes k})$, we define
\begin{equation*}
	\Phi_{f, k, m}(s) := P_{f, k, \leq m} s - T_{f, k} s \in \Gamma(X, L^{\otimes k}).
\end{equation*}
By the orthogonality relations in Theorem \ref{theorem: orthogonality-relations},
\begin{equation*}
	T_{f, k} s =\Pi_k(f\cdot s)= \Pi_k ( P_{f, k, \leq m} s) \in H^0(X, L^{\otimes k}).
\end{equation*}
Therefore, $\Phi_{f, k, m}(s)$ is orthogonal to the subspace $H^0(X, L^{\otimes k})$ of holomorphic sections.

By the Kodaira-H\"ormander $L^2$ estimates  in Theorem \ref{theorem: Hormander-estimate} (i.e., Theorem 2.2 in \cite{Cha2015}), there exists $k_0 \in \mathbb{Z}^+$ and a constant $C > 0$ (which are independent of $k, m, f, s$) such that, if $k \geq k_0$, then
\begin{align*}
	\lVert \Phi_{f, k, m}(s) \rVert & \leq C\cdot k^{-\frac{1}{2}}\cdot \lVert \bar{\partial}\Phi_{f, k,m}(s) \rVert\\
	&= C\cdot k^{-\frac{1}{2}} \left\lVert \overline{\partial} \left( P_{f, k, \leq m} s \right) \right\rVert\\
	&\leq C\cdot\tilde{C}\cdot k^{-\frac{m+1}{2}}\cdot \lVert s \rVert.
\end{align*}
In the second line, we have used the fact that $\bar{\partial}(T_{f,k}s)=0$, and in the last line, we have used Proposition \ref{proposition: estimate-partial-bar-D-f-action-s}. Finally, to account for the case $k < k_0$, we define $C_{f, m}$ as
\begin{equation*}
	C_{f, m} = \max \left\{ \lVert \Phi_{f, 1, m} \rVert_{op} \cdot 1^{\frac{m+1}{2}}, ..., \lVert \Phi_{f, k_0, m} \rVert_{op} \cdot k_0^{\frac{m+1}{2}}, C \cdot \tilde{C} \right\}.
\end{equation*}
\end{proof}

\begin{rmk}
In view of Proposition \ref{proposition: commutator-Bergman-projection-differential-operators}, Theorem \ref{theorem: main} provides a generalization of the result in Charles--Polterovich \cite{ChaPol2017}. 
\end{rmk}

\begin{rmk}\label{remark: module-version-of-BMS}
    Here we give an interpretation of Theorem \ref{theorem: main} by a comparison with the estimate of a composition of Berezin--Toeplitz operators in equation \eqref{equation: BT-star-product-asymptotic}. 
    On the one hand, for any $f\in C^\infty(X)$ and $s\in\mathcal{H}_k$, the output $P_{f,k,m}(s)$ depends only on the jets of $s$ at each point $x_0\in X$ (by locality of differential operators). 
    On the other hand, the operator $P_{f,k,m}|_{x_0}$ also depends on the jets of $f$ at $x_0$ by the Fedosov construction. Therefore, Theorem \ref{theorem: main} can be informally reformulated in the following way: 
    there exist bi-differential operators $\tilde{C}_l: C^\infty(X)\times\Gamma(X, L^{\otimes k})\rightarrow\Gamma(X, L^{\otimes k})$, such that 
    $$
    T_{f,k}(s)\sim\sum_{l\geq 1}\tilde{C}_l(f,s).
    $$
    We can hence regard Theorem \ref{theorem: main} as the ``module version'' of the asymptotic composition formula \eqref{equation: BT-star-product-asymptotic} from \cite{Guillemin95, BorMeiSch1994}. 
\end{rmk}



\subsection{The Fedosov-to-Berezin--Toeplitz degeneration}
\


In this subsection, we show that for the symbol of a  quantizable function, its associated Kostant-Souriau operators and Berezin-Toeplitz operator are identical. 
We first have the following proposition, which we call the {\em Fedosov-to-Berezin--Toeplitz degeneration}:
\begin{prop}
	Suppose $f$ is a formal quantizable  function on $X$. Then the Berezin--Toeplitz operator associated to its evaluation at 
	$\hbar=\tfrac{\sqrt{-1}}{k}$ is identical to the holomorphic differential operator defined via the Bargmann--Fock action. 
\end{prop}
\begin{proof}
	We  can decompose the formal function $f$ according to the power expansion of $\hbar$:
	$$
	f=f_0+\hbar\cdot f_1+\cdots +\hbar^l\cdot f_l,
	$$
	where each $f_i\in C^\infty(X)$. 
	The associated flat section $O_f$ are denoted by the stars in the following picture:
	\begin{equation}\label{diagram: formal-quantizable-function}
		\begin{sseq}[entrysize=1cm]{0...3}{0...3}
			\ssdrop{f_0}
			\ssmove{1}{0}
			\ssdrop{\hbar f_1}
			\ssmove 1 0
			\ssdrop{\hbar^2f_2}
			\ssmove 1 0
			\ssdrop{\cdots}
			\ssmoveto{0}{1}
			\ssdrop{*}
			\ssmove 1 0
			\ssdrop{*}
			\ssmove 1 0
			\ssdrop{*}
			\ssmove 1 0
			\ssdrop{\cdots}
			\ssmoveto{0}{2}
			\ssdrop{\cdots}
			\ssmove{1}{0}
			\ssdrop{\cdots}
			\ssmove{1}{0}
			\ssdrop{\cdots}
			\ssmove{1}{0}
			\ssdrop{\cdots}
		\end{sseq}
	\end{equation}
	The assumption that $f$ is formally quantizable is equivalent to saying that the flat section $O_f$ lives in a finite rectangle in the above picture. In particular, this finiteness condition guarantees that if we can take the evaluation 
	$\hbar=\tfrac{\sqrt{-1}}{k}$, there are no convergence issues. 
	
	For each $f_i$ and each weight $m\geq 1$, the orthogonality relation in Theorem \ref{theorem: orthogonality-relations} implies that for any holomorphic section $s\in\mathcal{H}_k$, we have
	\begin{equation}\label{equation: termwise-orthogonal}
		\sigma_{L^{\otimes k}} \left((O_{f_i})_m\circledast_k O_s\right)\in\mathcal{H}_k^{\perp}.
	\end{equation}
	
	On the one hand, by taking the sum of equation \eqref{equation: termwise-orthogonal} for all $i$ and  $m\geq 1$, we obtain that
	$$
	\sigma_{L^{\otimes k}} \left((O_{f}-f)\circledast_k O_s\right)\in\mathcal{H}_k^{\perp}.
	$$
	Here the Bargmann--Fock action $O_f\circledast_k O_s$ is well-defined since $f$ is formal quantizable. Also recall that in equation \eqref{equation: Wick-ordering-formula} which defines the level $k$ Bargmann--Fock action, we need to take the evaluation 
	$\hbar=\tfrac{\sqrt{-1}}{k}$.  
	
	On the other hand, we have seen in Proposition  \ref{proposition: Bargmann-Fock-differential-operator} that $\sigma_{L^{\otimes k}}\left(O_f\circledast_k O_s\right)\in\mathcal{H}_k$ since $s \mapsto \sigma_{L^{\otimes k}}(O_f\circledast_k O_s)$ is a holomorphic differential operator. Also, we obtain an orthogonal decomposition for $f\cdot s$:
	$$
	f\cdot s=-\sigma_{L^{\otimes k}}\left((O_f-f)\circledast_k s\right)+\sigma_{L^{\otimes k}}(O_f\circledast_k s).
	$$
	It follows from the definition of Berezin--Toeplitz operators that 
	$$
	T_{f,k}(s)=\Pi_k(f\cdot s)=\sigma_{L^{\otimes k}}(O_f\circledast_k s).
	$$
\end{proof}
\begin{rmk}
	We call this result the {\em Fedosov-to-Berezin--Toeplitz degeneration} for the following reason. Given a formal quantizable function $f=f_0+\hbar\cdot f_1+\hbar^2\cdot f_2+\cdots \hbar^m\cdot f_m$, on the one hand, the Berezin--Toeplitz operator $T_{f|_{\hbar=\sqrt{-1}/k},k}$ is defined using the terms in the first row of diagram \eqref{diagram: formal-quantizable-function}, i.e., multiplication by the function $f|_{\hbar=\sqrt{-1}/k}$ itself. On the other hand, the corresponding holomorphic differential operator is defined using the Bargmann--Fock action $s \mapsto \sigma_{L^{\otimes k}}(O_f \circledast_k O_s)$. However, all the terms in $O_f$ which are {\em not} in the first row will only contribute terms that are orthogonal to the subspace 
	$\mathcal{H}_k\subseteq \Gamma(X, L^{\otimes k})$. In other words, the holomorphic differential operator action (Fedosov approach) degenerates to the Berezin--Toeplitz operator corresponding to the first row 
	$T_{f_0+\sqrt{-1}/kf_1+(\sqrt{-1}/k)^2f_2+\cdots+(\sqrt{-1}/k)^mf_m}$. 
\end{rmk}


More generally, even if a level $k$ quantizable function is not globally the evaluation of the flat section associated to a formal quantizable function at 
$\hbar=\tfrac{\sqrt{-1}}{k}$, the statement still holds. We have the following theorem which tells us exactly when Berezin--Toeplitz operators act on the Hilbert spaces as genuine holomorphic differential operators.
\begin{thm}\label{theorem: Toeplitz-holomorphic-differential-operator}
	Suppose a smooth function $f$ is the symbol of a quantizable function of level $k$, i.e., $f=\sigma(\alpha)$ for a $D_{BT,k}$-flat section $\alpha\in\Gamma(X, \W_X \otimes \Sym \overline{T^*X})$. Then the corresponding Berezin--Toeplitz operator $T_{f,k}$ is equal to the holomorphic differential operator $P_{\alpha,k}$ on $\mathcal{H}_k$.
	
	Conversely, suppose a Berezin--Toeplitz operator $T_k:\mathcal{H}_k\rightarrow\mathcal{H}_k$ is a holomorphic differential operator, then it can be written as $T_k=T_{f,k}$, where $f$ is the symbol of a level $k$ quantizable function. 
\end{thm}
\begin{proof}
	For one direction, 
	by Proposition \ref{proposition: quantizable-function-locally-formal},  locally on an open set $U\subset X$, $\alpha$ must be the evaluation of the flat section associated to a formal quantizable function at 
	$\hbar=\tfrac{\sqrt{-1}}{k}$. Explicitly, there exists a formal quantizable function $f_0+\hbar\cdot f_1+\cdots+\hbar^l\cdot f_l\in C^\infty(U)[\hbar]$, such that 
	$$
	\alpha\vert_U=\left. \left(O_{f_0+\hbar\cdot f_1+\cdots+\hbar^l\cdot f_l}\right) \right\vert_{\hbar=\frac{\sqrt{-1}}{k}}.
	$$ 
	By taking symbols to both sides, we obtain that 
	$$
	f\vert_U= \left. \left(f_0+\hbar\cdot f_1+\cdots+\hbar^l\cdot f_l\right) \right\vert_{\hbar=\frac{\sqrt{-1}}{k}}.
	$$
	For any holomorphic sections $s,s'\in\mathcal{H}_k$, the proof of Theorem \ref{theorem: orthogonality-relations} implies that there exists a $2n-1$ form $\beta_U$ on $U$ such that
	\begin{equation}\label{equation: term-wise-exact}
		h( P_{\alpha,k}(s)-f\cdot s, s') \cdot\omega^n=\sum_{i=1}^l\sum_{m\geq 1} h\left( \left( \frac{\sqrt{-1}}{k} \right)^i \cdot\sigma_{L^{\otimes k}}\left((O_{f_i})_m\circledast_k s\right), s' \right)\cdot\omega^n=d\beta_U.
	\end{equation}
	Here $h$ denotes the hermitian metric on $L^{\otimes k}$. We claim that these anti-derivatives $\beta_U$ actually glue to a global one, which implies that the left hand side of equation \eqref{equation: term-wise-exact} is a globally exact $2n$-form on $X$. It follows that $P_{\alpha,k}(s)-f\cdot s\in\mathcal{H}_k^{\perp}$, and we obtain that
	$$
	P_{\alpha,k}(s)=T_{f,k}(s). 
	$$ 
	
	To prove this claim, first notice that on each open subset $U$ of $X$, different choices of local formal quantizable functions differ by a function in $(\hbar-\tfrac{\sqrt{-1}}{k})\cdot C^\infty(U)[\hbar]$, i.e.,  the ideal generated by $\hbar-\tfrac{\sqrt{-1}}{k}$. 
	A simple observation is that in the proof of Theorem \ref{theorem: orthogonality-relations},  the operations involved in finding the anti-derivatives preserves the ideal generated by $\hbar-\tfrac{\sqrt{-1}}{k}$. Thus we have 
	$$
	\beta_U-\beta_V\in(\hbar-\tfrac{\sqrt{-1}}{k})\cdot\A_U^{2n-1}[\hbar].
	$$
	Thus these local anti-derivatives coincide after taking the evaluation 
	$\hbar=\tfrac{\sqrt{-1}}{k}$ and glue together to a global one. 
	
	For the converse direction, notice that by the isomorphism in equation \eqref{equation: isomorphism-quantizable-function-differential-operator}, the differential operator corresponding to $T_k$ can be expressed as a Bargmann--Fock action $P_{\alpha,k}=\sigma_{L^{\otimes k}}\circ (\alpha\circledast_k-)$ for some section $\alpha$ of the Weyl bundle which is flat under $D_{BT,k}$. It follows that $T_k=P_{\alpha,k}=T_{\sigma(\alpha),k}$. 
\end{proof}
\begin{cor}\label{corollary: holomorphic-differential-operator-symbol}
	Suppose $\alpha_1$ and $\alpha_2$ are two level $k$ quantizable function which have the same symbol, i.e., $\sigma_{L^{\otimes k}}(\alpha_1)=\sigma_{L^{\otimes k}}(\alpha_2)=f$, then there corresponding differential operators are identical in $\End(\mathcal{H}_k)$. 
\end{cor}
\begin{proof}
	This follows immediately from Theorem \ref{theorem: Toeplitz-holomorphic-differential-operator}, since they are both $T_{f,k}$.
\end{proof}

\appendix

\section{Fedosov connections via quantization of $L_\infty$ structure}\label{appendix: Fedosov-connection-via-L-infinity-structure}

In this appendix, we briefly recall the explicit expression of the Fedosov connections $D_{BT}$ obtained via a quantum extension of Kapranov's $L_\infty$ structure on K\"ahler manifolds. For more details, we refer to \cite{ChaLeuLi2022b}.

\subsection{Kapranov's $L_\infty$ structure on K\"ahler manifolds}
\

The $L_\infty$ structure on K\"ahler manifolds was constructed by Kapranov in \cite{Kap1999}, which can be expressed in terms of a square zero differential on $\A_X^\bullet\otimes\W_X$ making it the corresponding Chevalley--Eilenberg complex. We will denote this differential by $D_K$, which is given explicitly as follows: first of all, there exist
$$
R_m^*\in\mathcal{A}^{0,1}_X(\Hom(T^*X,\Sym^m(T^*X))),\qquad m\geq 2
$$
such that their extensions $\tilde{R}^*_m$ to the holomorphic Weyl bundle $\W_X$ by derivation (with respect to the classical commutative product) satisfy
$$\left(\bar{\partial}+\sum_{m\geq 2}\tilde{R}_m^*\right)^2 = 0.$$

These $R_m^*$'s are defined as partial transposes of the higher covariant derivatives of the curvature tensor:
$$
R_2^*=\frac{1}{2}R_{i\bar{j}k}^p d\bar{z}^j\otimes (y^iy^k\otimes\partial_{z^p}),\qquad R_m^*=(\delta^{1,0})^{-1}\circ\nabla^{1,0}(R_{m-1}^*),
$$
where, by abuse of notations, we use $\nabla$ to denote the Levi-Civita connection on the (anti)holomorphic (co)tangent bundle of $X$, as well as their tensor products including the Weyl bundles. 
We can write these $R_m^*$'s locally as:
\begin{equation}\label{equation: terms-L-infty-structure}
	R_m^*=R_{i_1\cdots i_m,\bar{l}}^j d\bar{z}^l\otimes (y^{i_1}\cdots y^{i_m}\otimes \partial_{z^j}).
\end{equation}
Readers are referred to \cite[Section 2.5]{Kap1999} for a detailed exposition. The following opeator 
$$
D_K=\nabla-\delta^{1,0}+\sum_{m\geq 2}\tilde{R}_m^*
$$
squares to $0$ and defines a flat connection on the holomorphic Weyl bundle $\W_X$, which is compatible with the classical (commutative) product. The symbol map induces an isomorphism between (local) flat sections of $\W_X$ under $D_K$ and holomorphic functions:
$$
\sigma: \Gamma^{flat}(U, \W_X)\overset{\cong}{\rightarrow}\mathcal{O}_X(U).
$$
For a holomorphic function $f\in\mathcal{O}_X(U)$, the corresponding flat section is explicitly given by
\begin{equation}\label{equation: flat-section-holomorphic-function}
O_f=\sum_{m\geq 0}(\tilde{\nabla}^{1,0})^m(f).
\end{equation}
In particular, $O_f\in\Gamma(U, \W_X)$ and is of polarized weight $0$.

\subsection{Quantizing $D_K$ to Fedosov connections}\label{subsection: Fedosov-quantization}
\

To quantize the flat connection $D_K$, roughly speaking, we turn the terms $\tilde{R}_m^*$ to brackets with respect to the fiberwise Wick product. In particular, we have to replace $\W_X$ by the formal complexified Weyl bundle $\W_{X,\C}[[\hbar]]$. 

More explicitly, we consider the $\A_X^{\bullet}$-linear operator
$$L:\A_X^\bullet\left(\widehat{\Sym}(T^*X)\otimes TX\right)\rightarrow \A_X^\bullet\left(\widehat{\Sym}(T^*X)\otimes \overline{T^*X}\right)$$
given by ``lifting the last subscript'' using the K\"ahler form $\omega$. Then we define
$$
I_m:=L(R_m^*)=R_{i_1\cdots i_m,\bar{l}}^j\omega_{j\bar{k}}d\bar{z}^l\otimes (y^{i_1}\cdots y^{i_m}\bar{y}^{k})\in\mathcal{A}_X^{0,1}(\mathcal{W}_{X,\mathbb{C}}). 
$$
Let $I:=\sum_{m\geq 2}I_m$. We define a connection on $X$ by
$$
D_F:=\nabla-\delta+\frac{1}{\hbar}[I,-]_\star.
$$
It was shown in \cite{ChaLeuLi2022b} that $D_F^2=0$, and thus defines a Fedosov connection on $\W_{X,\C}$.
We say that $D_F$ is a {\em quantum extension} of $D_K$ because $D_F|_{\W_X}=D_K$.

The Fedosov connection $D_F$ induces a deformation quantization of $X$, which is of Wick type. It was shown by Karabegov in \cite{Kar1996} that Wick type deformation quantizations on a K\"ahler manifold $X$ are classified by formal closed 2-forms on $X$ of type $(1,1)$, known as the {\em Karabegov form}. For instance, the Karabegov form of the Berezin--Toeplitz deformation quantization is given by
$$
-\frac{1}{\hbar}\omega+\text{Ric}_X. 
$$

It turns out that for all Wick type deformation quantization of $X$, there exists a quantum extension of $D_K$ inducing this star product. For every formal closed $2$-form $\alpha\in\A_X^{1,1}[[\hbar]]$ of type $(1,1)$, there exists a Fedosov connection $D_{F,\alpha}$ such that 
\begin{enumerate}
	\item $D_{F,\alpha}^2=0$;
	\item $D_{F,\alpha}|_{\W_X}=D_K$;
	\item The Karabegov form of the associated Wick type deformation quantization on $X$ is $-\frac{1}{\hbar}\omega+\alpha$. 
\end{enumerate}
This connection is of the form 
$$
D_{F,\alpha}=\nabla-\delta+\frac{1}{\hbar}[I+J_\alpha,-]_\star.
$$
Here the terms $J_\alpha$ are defined as follows: the $\partial\bar{\partial}$-lemma guarantees the local existence of a formal function $g$ such that $\alpha=\bar{\partial}{\partial}(g)$. We define $J_\alpha\in \mathcal{A}_X^{0,1}(\mathcal{W}_X)$ by
\begin{equation*}\label{equation: I-alpha-formula}
	J_\alpha:=\hbar\cdot\sum_{k\geq 1}\left((\delta^{1,0})^{-1}\circ\nabla^{1,0}\right)^k(\bar{\partial}g).
\end{equation*} 
It is easy to see that $J_\alpha$ is independent of the choice of the potential $g$. 
For all these $\alpha$'s, the connections $D_{F,\alpha}$ behave nicely because of the explicit forms of the terms $I$ and $J_{\alpha}$'s. Specifically, $J_{\alpha}$ only contains polynomials in $\W_X$, and $I\in\W_X\otimes \overline{TX}^*$, i.e., it is degree $1$ in $\overline{\W_X}$.

For the Berezin--Toeplitz star product $\star_{BT}$, we denote the corresponding Fedosov connection by 
\begin{equation}\label{equation: Fedosov-connection-Toeplitz-quantization}
	D_{BT}:=\nabla-\delta+\frac{1}{\hbar}[I_{BT},-]_\star=\nabla+\frac{1}{\hbar}[\gamma_{BT},-]_\star.
\end{equation}
Here the second equality follows from the fact that the operator $\delta$ can also be written as a bracket. A simple counting of polynomial degrees of $\hbar$ and $\bar{y}^j$'s gives the following
\begin{lem}
	The term $\frac{1}{\hbar}\cdot[I_{BT},-]_\star$ in the Fedosov connection $D_{BT}$ has weight $0$ with respect to the polarized weight defined in Definition \ref{definition: polarized-weight}. 	
\end{lem}

\subsection{Non-formal Fedosov connection}\label{subsection: non-formal-Fedosov-connection}
\

For any quantum level $k>0$, we can take the evaluation $\hbar=\sqrt{-1}/k$ to obtain a non-formal Fedosov connection 
\begin{equation}\label{equation: level-k-Fedosov-connection}
D_{BT,k}=\nabla-\delta+\frac{k}{\sqrt{-1}}\cdot[I_{BT}|_{\hbar=\sqrt{-1}/k},-]_{\star_k}.
\end{equation}
Here we also replace the formal Wick product by the level $k$ Wick product. 
\begin{lem}
 The connection $D_{BT,k}$ is a well-defined flat connection on the sub-bundle $\W_X\otimes \Sym\overline{T^*X}$ of the complexified Weyl bundle.
\end{lem}
\begin{proof}
The flatness $D_{BT,k}^2=0$ follows from $D_{BT}^2=0$ by taking the evaluation $\hbar=\sqrt{-1}/k$. We only need to show that it is well-defined. The main issue here is to avoid an infinite power series of $1/k$: when restricted to the sub-bundle $\W_X\otimes\Sym\overline{T^*X}$, this is indeed the case, since $I_{BT}$ is a power series only in $T^*X$ which contracts with $\overline{T^*X}$ in the level $k$ Wick product $\star_k$. 
\end{proof}

\bibliographystyle{amsplain}
\bibliography{References}

\end{document}